\documentclass[a4paper,reqno,11pt]{amsart}

\newcommand{\myfont}{\sffamily}

\usepackage[normalem]{ulem}

\newtheoremstyle{mythmstyle}
  {\topsep}
  {\topsep}
  {\itshape}
  {}
  {\bfseries \myfont}
  {.}
  {.5em}
  {}

\newtheoremstyle{mydefstyle}
  {\topsep}
  {\topsep}
  {\normalfont}
  {}
  {\bfseries \myfont}
  {.}
  {.5em}
  {}

\swapnumbers                    

\theoremstyle{mythmstyle}       

\let\MakeUppercase\relax

\expandafter\let\expandafter\oldproof\csname\string\proof\endcsname
\let\oldendproof\endproof
\renewenvironment{proof}[1][\bfseries\myfont\proofname]{%
  \oldproof[\bfseries \myfont #1]%
}{\oldendproof}

\makeatletter
\newtoks\thm@headfont  \thm@headfont{\bfseries \myfont}

\renewcommand\section{\@startsection{section}{1}%
  \z@{.7\linespacing\@plus\linespacing}{.5\linespacing}%
  {\Large\myfont\bfseries}}
\renewcommand\subsection{\@startsection{subsection}{2}%
  \z@{-.5\linespacing\@plus-.7\linespacing}{.5\linespacing}%
  {\large\myfont\bfseries}}
\renewcommand\subsubsection{\@startsection{subsubsection}{3}%
  \z@{.5\linespacing\@plus.7\linespacing}{-.5em}%
  {\myfont\bfseries}}
\renewenvironment{abstract}{%
  \ifx\maketitle\relax
    \ClassWarning{\@classname}{Abstract should precede
      \protect\maketitle\space in AMS document classes; reported}%
  \fi
  \global\setbox\abstractbox=\vtop \bgroup
    \normalfont\Small
    \list{}{\labelwidth\z@
      \leftmargin3pc \rightmargin\leftmargin
      \listparindent\normalparindent \itemindent\z@
      \parsep\z@ \@plus\p@
      
    }%
    \item[\hskip\labelsep
      \myfont\bfseries
    \abstractname.]%
}{%
  \endlist\egroup
  \ifx\@setabstract\relax \@setabstracta \fi
}
\renewcommand\contentsnamefont{\myfont\bfseries}
\renewcommand\@starttoc[2]{\begingroup
  \setTrue{#1}%
  \par\removelastskip\vskip\z@skip
  \@startsection{}\@M\z@{\linespacing\@plus\linespacing}%
    {.5\linespacing}{
      \contentsnamefont}{#2}%
  \ifx\contentsname#2%
  \else \addcontentsline{toc}{section}{#2}\fi
  \makeatletter
  \@input{\jobname.#1}%
  \if@filesw
    \@xp\newwrite\csname tf@#1\endcsname
    \immediate\@xp\openout\csname tf@#1\endcsname \jobname.#1\relax
  \fi
  \global\@nobreakfalse \endgroup
  \addvspace{32\p@\@plus14\p@}%
  \let\tableofcontents\relax
}
\renewcommand\@settitle{\begin{center}%
  \baselineskip14\p@\relax
    \LARGE
    \bfseries
    \myfont
  \@title
  \end{center}%
}
\renewcommand\@setauthors{%
  \begingroup
  \def\thanks{\protect\thanks@warning}%
  \trivlist
  \centering\footnotesize \@topsep30\p@\relax
  \advance\@topsep by -\baselineskip
  \item\relax
  \author@andify\authors
  \def\\{\protect\linebreak}%
  \large
  \myfont\bfseries\authors
  \ifx\@empty\contribs
  \else
    ,\penalty-3 \space \@setcontribs
    \@closetoccontribs
  \fi
  \endtrivlist
  \normalfont\myfont\@setaddresses
  \endgroup
}
\renewcommand\@setaddresses{\par
  \nobreak \begingroup
\footnotesize
  \def\author##1{\nobreak\addvspace\bigskipamount}%
  \def\\{\unskip, \ignorespaces}%
  \interlinepenalty\@M
  \def\address##1##2{\begingroup
    \par\addvspace\bigskipamount\indent
    \@ifnotempty{##1}{(\ignorespaces##1\unskip) }%
    {
      \ignorespaces##2}\par\endgroup}%
  \def\curraddr##1##2{\begingroup
    \@ifnotempty{##2}{\nobreak\indent\curraddrname
      \@ifnotempty{##1}{, \ignorespaces##1\unskip}\/:\space
      ##2\par}\endgroup}%
  \def\email##1##2{\begingroup
    \@ifnotempty{##2}{\nobreak\indent\emailaddrname
      \@ifnotempty{##1}{, \ignorespaces##1\unskip}\/:\space
      \ttfamily##2\par}\endgroup}%
  \def\urladdr##1##2{\begingroup
    \def~{\char`\~}%
    \@ifnotempty{##2}{\nobreak\indent\urladdrname
      \@ifnotempty{##1}{, \ignorespaces##1\unskip}\/:\space
      \ttfamily##2\par}\endgroup}%
  \addresses
  \endgroup
}
\renewcommand\enddoc@text{\ifx\@empty\@translators \else\@settranslators\fi
}
\renewcommand\@secnumfont{\myfont\bfseries} 

\renewcommand\maketitle{\par
  \@topnum\z@ 
  \@setcopyright
  \thispagestyle{firstpage}
  \ifx\@empty\shortauthors \let\shortauthors\shorttitle
  \else \andify\shortauthors
  \fi
  \@maketitle@hook
  \begingroup
  \@maketitle
  \toks@\@xp{\shortauthors}\@temptokena\@xp{\shorttitle}%
  \toks4{\def\\{ \ignorespaces}}
  \edef\@tempa{%
    \@nx\markboth{\the\toks4
      \@nx\MakeUppercase{\the\toks@}}{\the\@temptokena}}%
  \@tempa
  \endgroup
  \c@footnote\z@
  \@cleartopmattertags
}
\def\@captionheadfont{\myfont\bfseries} 

\makeatother

  \usepackage{datetime}    
  
  \yyyymmdddate   
  \AddToHook{shipout/background}[mypic]
  {
    \put(2.5cm,-28.5cm)%
    {\footnotesize \today, \currenttime, \emph{File:} \texttt{\jobname.tex}}

  }
\usepackage[a4paper]{geometry}
\usepackage{enumitem}
\setlist[enumerate,1]{label={\upshape(\alph*)}}
\setlist[enumerate]{itemsep=1.ex,topsep=0.8ex,leftmargin=5ex,labelwidth=5ex}

\usepackage{lmodern}

\usepackage{amssymb,amsmath,amsthm}
\usepackage{mathtools}
\usepackage{mathrsfs}
\usepackage{bbm}

\usepackage[font=small]{caption}

\usepackage{color}
\definecolor{darkred}{rgb}{0.2,0.2,0.6}
\definecolor{darkblue}{rgb}{0.2,0.2,0.6}
\usepackage[colorlinks=true,linkcolor=blue, citecolor=darkred]{hyperref}
\usepackage{graphicx}

\usepackage{accents}  
\newcommand\sfA{{\mathsf A}}

 \newcommand\sfH{{\mathsf H}}
\newcommand\sfI{{\mathsf I}} \newcommand\sfJ{{\mathsf J}}
\newcommand\sfK{{\mathsf K}} \newcommand\sfL{{\mathsf L}}

 \newcommand\sfR{{\mathsf R}}
 
\newcommand\sfU{{\mathsf U}}

\newcommand\fra{\mathfrak a}

\newcommand\dC{{\mathbb C}}

 \newcommand\dN{{\mathbb N}}
 
 \newcommand\dR{{\mathbb R}}

 \newcommand\dZ{{\mathbb Z}}

 \newcommand\cB{{\mathcal B}}

\newcommand\ii{{\mathsf{i}}}  

\newcommand\eps{\varepsilon}
\renewcommand\phi{\varphi}

\DeclareMathOperator{\dist}{dist}
\DeclareMathOperator{\dom}{dom}
\DeclareMathOperator{\ran}{ran}

\DeclareMathOperator{\supp}{supp}

\DeclareMathOperator{\dd}    {d\!}

\newcommand{\restr}[1]{{\restriction}_{#1}} 
\newcommand{\conj}[1]{\overline {#1}}  
\newcommand{\defeq}{\vcentcolon=}

\newcommand{\bd}  {\partial}          
\newcommand{\clo}[2][]{\overline{{#2}}^{#1}} 

\DeclarePairedDelimiter\abs{\lvert}{\rvert}%
\DeclarePairedDelimiter\norm{\lVert}{\rVert}%
\DeclarePairedDelimiter\myiprod{\langle}{\rangle}%
\newcommand{\iprod}[3][]{\myiprod[#1]{#2, #3}}%

\newcommand{\Lsqr}{L_2}                  %
\newcommand{\Linfty}{L_\infty}                  %
\newcommand{\Ci}{C^\infty}  
\newcommand{\Cci}{C^\infty_{\mathrm c}}  
\newcommand{\Sob}{H}        
\newcommand{\Sobn}{\ring H}  
\newcommand\ov{\overline}
\newcommand\wt{\widetilde}

\newcounter{counter_a}

\newcommand{\itemref}[1]{\noindent \ref{#1}}

\newcommand{\eg}{\textit{e.g.\ }}
\newcommand{\ie}{\textit{i.e.\ }}
\newcommand{\cf}{\textit{cf.\ }}

\newcommand{\quadtext}[1]{\quad\text{#1}\quad}
\newcommand{\qquadtext}[1]{\qquad\text{#1}\qquad}

\numberwithin{figure}{section}
\numberwithin{equation}{section}

\swapnumbers  
\newcounter{intro}

\newtheorem{maintheorem}[intro]{Theorem}

\swapnumbers 
\newtheorem{theorem}{Theorem}[section]
\newtheorem*{theorem*}{Theorem}
\newtheorem{proposition}[theorem]{Proposition}
\newtheorem{lemma}[theorem]{Lemma}
\newtheorem{corollary}[theorem]{Corollary}

\theoremstyle{definition}       
\newtheorem{definition}[theorem]{Definition}

\newtheorem{example}[theorem]{Example}

\newtheorem{remark}[theorem]{Remark}

\newtheorem*{remark*}{Remark}

\newtheorem{hypothesis}[theorem]{Hypothesis}

\newcommand{\R}{\mathbb{R}}

\begin{document}

\title[]{Convergence of Schr\"odinger operators on domains with scaled
  resonant potentials}

\author[V. Lotoreichik]{Vladimir Lotoreichik} %
\address[V.~Lotoreichik]{Department of Theoretical Physics, Nuclear
  Physics Institute, \newline
  Czech Academy of Sciences, Hlavn\'{i} 130, 25068
  \v Re\v z, Czech Republic} %
\email{lotoreichik@ujf.cas.cz} %
\author[O. Post]{Olaf Post} %
\address[O.~Post]{Fachbereich 4 – Mathematik, Universit\"{a}t Trier,
  54286 Trier, Germany} %
\email{olaf.post@uni-trier.de}



\begin{abstract}
  We consider Schr\"odinger operators on a bounded, smooth domain of
  dimension $d \ge 2$ with Dirichlet boundary conditions and a
  properly scaled potential, which depends only on the distance to the
  boundary of the domain. Our aim is to analyse the convergence of
  these operators as the scaling parameter tends to zero. If the
  scaled potential is resonant, the limit in strong resolvent sense is
  a Robin Laplacian with boundary coefficient expressed in terms of
  the mean curvature of the boundary.  A counterexample shows that
  norm resolvent convergence cannot hold in general in this
  setting. If the scaled potential is non-resonant and satisfies an explicit assumption on the smallness of the negative part, the limit in strong resolvent sense is the Dirichlet
  Laplacian. We conjecture that we can drop this additional assumption in the non-resonant case.
\end{abstract}
\maketitle

%
\section{Introduction}
\label{sec:intro}
%

\subsection{Background and motivation}
Convergence of Schr\"odinger operators with scaled potentials is a
classical topic in mathematical physics. A collection of results in
this area can be found in the monograph~\cite{AGHH}.  In
many cases, it is observed that one obtains in the strong or norm
resolvent limit an operator with an interaction supported on a set of
measure zero, being defined via a boundary condition. The structure of
the limiting operator often drastically depends on whether the scaled
potential is resonant or not.

A typical example, where the limit depends on whether the scaled
potential is resonant, is the approximation of a Schr\"odinger
operator with a point $\delta$-interaction in three dimensions by a
family of Schr\"odinger operators with scaled regular
potentials. First partial results on this approximation are obtained
in~\cite{AFH79, F72, N77}.  To the best of our knowledge
approximations of $\delta$-interactions were first addressed in full
detail by Albeverio and H{\o}egh-Krohn in~\cite{AH81}, who proved the
convergence in strong resolvent sense. The limit has a non-trivial
point interaction if the scaled potential is resonant, otherwise the
limit is the free Laplacian. Shortly thereafter, it was shown	
in~\cite{AGH82} that under rather non-restrictive additional assumptions \emph{norm} instead of \emph{strong} resolvent
convergence holds in this approximation.  By restricting to the class of radially
symmetric potentials one obtains upon separation of variables a model
problem on the half-line, which was considered separately by \v{S}eba
in~\cite{seba:85}.  He proved that a family of half-line Schr\"odinger
operators with Dirichlet boundary conditions and locally scaled
potentials converges in norm resolvent sense to the Neumann (or Robin)
Laplacian on the half-line if the scaled potential is resonant and to
the Dirichlet Laplacian otherwise.  Further refinements and extensions of these
results can be found in~\cite{DM16, SLS21}.

Typically, in the analysis of such a convergence, the integral kernel
of the resolvent of the unperturbed operator is used and the resolvent
identity plays a significant role, even though, in certain
approximation problems with a non-explicit integral kernel it suffices
to know only its singular part; \cf the recent analysis~\cite{NS24} of
approximation of point interactions on bounded domains.  Our main
motivation is to develop an approach to this class of problems, which
does not use the integral kernels of the resolvent, and where the
analysis is merely performed on the level of quadratic forms.  The
advantage of our method here is that it can be efficiently applied to
settings, in which the integral kernel of the resolvent of the
unperturbed operator is not given in an explicit form. Our final goal
and main motivation is to analyse the convergence of Schr\"odinger
operators on bounded smooth domains with Dirichlet boundary conditions
and scaled potentials, which depend only on the distance to the
boundary.

We remark that a similar phenomenon, where the limit qualitatively
depends on whether the scaled potential is resonant or not, was
also observed in full
 line Schr\"odinger operators with scaled potentials~\cite{ACF07, CE07} motivated by the applications to curved quantum waveguides. 
An extension of these results to scaled potentials with zero mean was later performed in a series of
papers~\cite{GM09, GH13,G22}, where the first two papers treat the
one-dimensional case, while the last one deals with the two-dimensional
case. In the latter papers, the scaled potentials converge to $\delta'$ in the sense of distributions. 
Thus, the resulting limit may then be interpreted as a rigorous definition of  a Schr\"odinger operator with a $\delta'$-potential.

Potentials, dependent only on the distance to a hypersurface, are also
used in the approximation of Schr\"odinger operators with surface
$\delta$-interactions~\cite{BEHL17, BEHL20, exner-ichinose:01} and
Dirac operators with $\delta$-shell interactions~\cite{BHS23, CLMT23,
  MP18}.  In these settings the effect of resonant potentials does not
occur and the choice of the potential merely manifests in the values
of the parameters characterising the limiting operator. In a certain
sense the limit ``continuously'' depends on the approximating
potential. Another important difference is that in those settings the
convergence typically holds in norm resolvent sense, while in the
setting considered in the present paper, in general, only strong
resolvent convergence can be proved, as a counterexample shows.

The proof of our main result for resonant potentials relies on the
construction of a suitable identification operator between the form
domains of the limiting operator and of the operators with scaled
resonant potentials. The key idea is to employ multiplication with the scaled resonant solution as cut-off function in such identification
operators.  We analyse the non-resonant case only partially and use a
similar method, in which we employ instead the derivative of the
non-resonant solution in the construction of identification operators.

\subsection{Resonant potentials in  one dimension}
We use the definition of resonant potentials borrowed
from~\cite{seba:85}. This class will be used throughout the whole paper.
\begin{definition}
  \label{def:resonant}
  The real-valued potential $V\in \Cci(\clo \dR_+)$ is called
  \emph{resonant} if the initial-value problem
  \begin{equation}
    \label{eq:IVP}
    \begin{cases}
      -\psi'' + V\psi = 0,\qquad\text{on}\,\,\dR_+,\\
      \psi(0) = 0,
    \end{cases}
  \end{equation}
  has a bounded non-trivial solution $\psi_0\in \Ci(\clo \dR_+)$
  (called \emph{resonant solution}). For the sake of convenience, we
  assume that $\supp V \subset[0,a]$ with some $a > 0$ and normalise
  the solution $\psi_0$ so that $\psi_0(t) = 1$ for all $t > a$.
\end{definition}

\begin{remark}
  \label{rem:resonant}
  Several observations on resonant potentials are in order.
  \begin{enumerate}
  \item
    \label{resonant.a}
    It is not hard to see that $V\in \Cci(\clo \dR_+)$ is resonant if
    and only if the Schr\"odinger operator with potential $V$ on the
    interval $(0,a)$ with Dirichlet boundary condition at $t = 0$ and
    Neumann boundary condition at $t = a$ has eigenvalue zero.
    Indeed, the continuous extension of the corresponding
    eigenfunction by a constant for $t > a$ gives a bounded solution
    of~\eqref{eq:IVP}. Conversely, the restriction of a bounded
    solution $\psi_0$ to $(0,a)$ is in the kernel of the
    aforementioned Schr\"odinger operator.

  \item
    \label{resonant.b}
    The potential $V\in \Cci(\clo \dR_+)$ not satisfying
    Definition~\ref{def:resonant} will be called
    \emph{non-resonant}. By the observation in (a) of this remark, we
    immediately see that any non-negative potential
    $V\in \Cci(\clo \dR_+)$ is non-resonant. For non-resonant
    potentials, we normalise the solution $\psi_0$ of the initial
    value problem~\eqref{eq:IVP} so that $\psi_0'(t) = 1$ for all
    $t > a$. We call such a solution $\psi_0$ \emph{non-resonant}.

  \item
    \label{resonant.c}
    If the potential $V$ is non-positive, then there exists a sequence
    $0 < \alpha_1 < \alpha_2 < \dots < \alpha_n < \dots$ with
    $\alpha_n\to\infty$ such that the potential $\alpha V$ with
    $\alpha > 0$ is resonant if and only if
    $\alpha\in\{\alpha_1,\alpha_2,\dots\}$.  The assumption on the
    smoothness of the potential $V$ is imposed for technical reasons
    only (\eg when using $V'$ in Lemma~\ref{lem:form-diff.a'}).  In
    particular, for the characteristic function $\chi_{(0,1)}$ of the
    interval $(0,1)$ the potential $V = -\alpha\chi_{(0,1)}$ is
    resonant if and only if $\alpha = (n+1/2)^2\pi^2$ for some
    $n\in\dN_0$ (\cf \cite{seba:85}).
  \end{enumerate}
\end{remark}
\begin{figure}[h]
  \includegraphics[width=6cm]{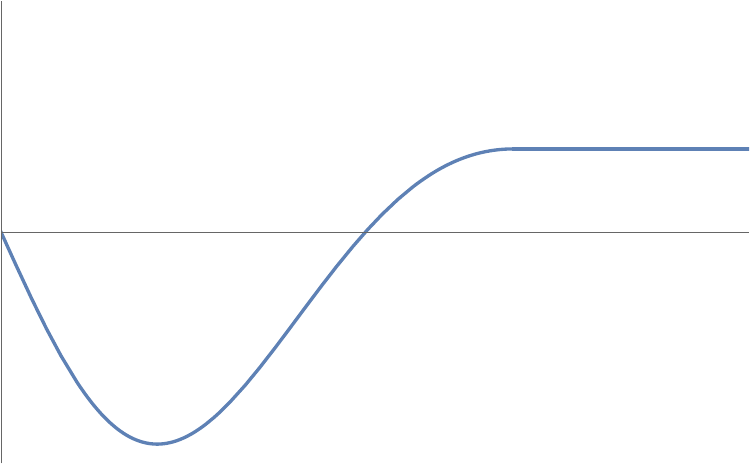} %
  \hspace{12ex} %
  \includegraphics[width=6cm]{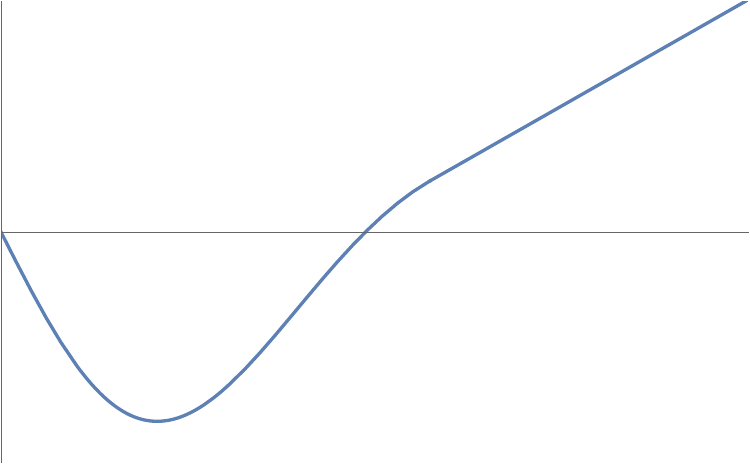}
  \caption{The solutions $\psi_0$ of~\eqref{eq:IVP} for resonant and
    non-resonant potentials on the left and on the right,
    respectively.}
\end{figure}
It was established in~\cite{seba:85} that the family of self-adjoint
Schr\"odinger operators in $\Lsqr(\dR_+)$ with Dirichlet boundary
conditions
\begin{equation*}
  \sfH_\eps \psi
  \defeq -\psi''
  + \frac{1}{\eps^2}V\Bigl(\frac{\cdot}{\eps}\Bigr)\psi,
  \qquad %
  \dom\sfH_\eps \defeq  \Sobn^1(\dR_+)\cap \Sob^2(\dR_+)
\end{equation*}
converges in norm resolvent sense (as $\eps\to0)$ to the self-adjoint
one-dimensional Laplace operator in $\Lsqr(\dR_+)$ with Neumann
boundary conditions
\begin{equation*}
  \sfH \psi \defeq -\psi'',
  \qquad \dom\sfH \defeq \big\{\psi\in \Sob^2(\dR_+)\colon \psi'(0) = 0\big\},
\end{equation*}
provided the potential $V$ is resonant.  However, if the potential $V$
is non-resonant the family of the operators $\sfH_\eps$ converges in
norm resolvent sense (as $\eps \to 0$) to the Dirichlet Laplacian
on the half-line
\begin{equation*}
  \sfH_0\psi \defeq -\psi'',
  \qquad \dom \sfH_0 \defeq  \Sobn^1(\dR_+)\cap \Sob^2(\dR_+).
\end{equation*}
There are two main features in this approximation problem.  Using
appropriate test functions of the form
$\psi_\eps(t)=\eps^{-1/2}\psi(t/\eps)$, it can be checked that the
operators in the approximating family are in general not uniformly
bounded from below.  Second, in the case of resonant potentials the
form domain $\Sobn^1(\dR_+)$ of the approximating operators is a
proper subspace of the form domain $\Sob^1(\dR_+)$ of the limit. These
two features make the approximation problem difficult to treat with
standard techniques based on comparison of quadratic forms;
\cf~\cite[Theorem~VI.3.6]{kato:66}.

Our aim in the present paper is to address a multi-dimensional
counterpart of these convergence results, which we will describe below
in detail. We remark that the result of \v{S}eba in~\cite{seba:85} is
more general than we stated here. He has also shown how to approximate the
one-dimensional Laplacian on the half-line with Robin boundary
conditions by means of replacing $1/\eps^2$ by $(1+\beta\eps)/\eps^2$
in the operator family $\sfH_\eps$\footnote{We remark that in the final formula in~\cite{seba:85} for the boundary parameter in the limiting operator the parameter $\beta$ 
appears in the denominator and formally the case $\beta = 0$ is excluded, but careful inspection of the proof shows that it applies also to the case $\beta = 0$ leading to Neumann boundary conditions in the limit.}. 
We will not address this more
general case in our analysis of the multi-dimensional problem, as we
already obtain a Robin-type boundary conditions by geometry.

\subsection{Main results}
Let $\Omega\subset\dR^d$, $d\ge 2$, be a bounded domain with
$\Ci$-smooth connected boundary $\Sigma \defeq \bd\Omega$.  The
$\Ci$-smoothness of the boundary and boundedness of the domain are
assumed for convenience and most of our analysis extends to
$C^2$-smooth domains with not necessarily compact boundaries under
additional uniformity assumptions on the boundary (\eg we have to
require $\rho>0$ in~\eqref{eq:mu}).  We denote by $\nu$ the outer unit
normal vector to $\Omega$.  The differential of the Gauss map
$\Sigma\ni s\mapsto \nu(s)$
\begin{equation}
  \label{eq:shape}
  \sfL_s \defeq \dd \nu(s)\colon T_s\Sigma \to T_s\Sigma
\end{equation}
is called the \emph{shape operator}; here $T_s\Sigma$ is the tangent
space for $\Sigma$ at $s\in\Sigma$.  The eigenvalues
$\kappa_1(s),\kappa_2(s),\dots,\kappa_{d-1}(s)\colon
\bd\Omega\to\dR$ of $\sfL_s$ are called the \emph{principal
  curvatures} of $\bd\Omega$. In our sign convention for the shape operator, all the
principal curvatures are non-negative if and only if $\Omega$ is
convex.  The mean curvature of $\Sigma$ at $s\in\Sigma$ is defined as
usual by
\begin{equation}
  \label{eq:meancurv}
  H(s) \defeq \frac{1}{d-1}\sum_{j=1}^{d-1}\kappa_{j}(s).
\end{equation}
Under our regularity assumptions on $\Omega$ the mean curvature $H$ is
a $\Ci$-smooth function on $\Sigma$. Note that the mean curvature
is non-negative for $\Omega$ being convex, while the converse is only
true in two dimensions.

We adopt the notation $\Sob^k(\Omega)$ for the $\Lsqr$-based Sobolev
space on $\Omega$ of order $k\in\dN$ and $\Sob^s(\bd\Omega)$ for the
$\Lsqr$-based Sobolev space on the boundary $\bd\Omega$ of order
$s\in\dR_+$. We denote by $\Sobn^k(\Omega)$ the closure of $C^\infty_{\rm c}(\Omega)$ in the Sobolev space $\Sob^k(\Omega)$ with respect to its norm.  
For a function $u\in \Sob^2(\Omega)$, we use the
notation $u\restr\Sigma \in \Sob^{3/2}(\Sigma)$ for its trace on the
boundary and $\partial_\nu u \restr \Sigma\in \Sob^{1/2}(\Sigma)$ for
the normal derivative corresponding to the normal pointing outwards of
$\Omega$.  The role of the one-dimensional Neumann Laplacian $\sfH$
from the previous subsection is played now by the self-adjoint Robin
Laplacian in $\Lsqr(\Omega)$
\begin{equation}
  \label{eq:A}
  \sfA u \defeq -\Delta u,\qquad
  \dom\sfA :=
  \Bigr\{u\in \Sob^2(\Omega)\colon \partial_\nu u\restr\Sigma
    + \frac{d-1}{2}H u\restr\Sigma = 0\Bigr\},
\end{equation}
where the Robin coefficient is expressed in terms of the mean
curvature. At the points where the mean curvature vanishes, we recover
locally Neumann boundary conditions. The role of the family of
one-dimensional Schr\"odinger operators $\sfH_\eps$ is played by the
self-adjoint Schr\"odinger operator in $\Lsqr(\Omega)$ defined for
$\eps > 0$ by
\begin{equation}
  \label{eq:Aeps}
  \sfA_\eps u
  \defeq -\Delta u
  + V_\eps u,
  \quad \dom \sfA_\eps := \Sobn^1(\Omega)\cap \Sob^2(\Omega),
  \quad\quadtext{where}
  V_\eps := \frac{1}{\eps^2}V \Bigl(\frac{\dist(\cdot\,,\Sigma)}{\eps}\Bigr)
\end{equation}
and where $V\in \Cci(\clo \dR_+)$.  Our first main result concerns the
class of resonant potential.
\begin{maintheorem}[the resonant case]
  \label{thm:main}
  Assume that $V\in \Cci(\clo \dR_+)$ is resonant in the
  sense of Definition~\ref{def:resonant}. Then the family of scaled
  Schr\"odinger operators $\sfA_\eps$ converges to the Robin Laplacian
  $\sfA$
  in strong resolvent sense as $\eps\to 0$.
\end{maintheorem}
In the proof of this result we rely on a convenient representation of
the sesquilinear form for the resolvent difference of the operators
$\sfA$ and $\sfA_\eps$ in terms of the resonant solution $\psi_0$. The
technique shares common ideas with the abstract approach for proving
norm resolvent convergence developed by the second-named author; see
the monograph~\cite{post:12} and the references therein.  Since the
form domains of $\sfA$ and $\sfA_\eps$ are different but the Hilbert
spaces are the same, one only needs identification operators mapping
from one form domain into the other.  Thus, the analysis boils down to
find a suitable identification operator, which maps a function
$\Sob^1(\Omega)$ into a function $\Sobn^1(\Omega)$.  In the
construction of this operator we use the resonant solution $\psi_0$
of~\eqref{eq:IVP}.  In Section~\ref{sec:counterex}, we construct a
counterexample, which shows that the operators $\sfA_\eps$ do not
converge to $\sfA$ in \emph{norm} resolvent sense.  This
counterexample relies on the analysis of the disk, where separation of
variables is available.  In particular, we cannot expect in general
that norm resolvent convergence holds in Theorem~\ref{thm:main}.
\begin{remark}[appearance of mean curvature terms in related problems]
  \indent
  \begin{enumerate}
  \item Note that the mean curvature term arises also in the large
    coupling asymptotics of the Robin Laplacian with a negative
    boundary parameter~\cite{EMP14, KP17, PP15, PP16} and for the
    Robin Laplacian on a shell in the small thickness
    limit~\cite{KRRS18}.
  \item The appearance of the curvature term in the boundary
    conditions of the limiting operator for scaled resonant potentials
    was also observed in~\cite{G22} in two dimensions for a different
    approximation problem, in which the limit has transmission
    $\delta^\prime$-type boundary conditions. The result there is
    proved by a different technique and is stated in terms of
    convergence of eigenvalues and weak convergence of eigenfunctions.
  \end{enumerate}
\end{remark}

The role of the one-dimensional Dirichlet Laplacian $\sfH_0$ from the
previous subsection is played now by the self-adjoint Dirichlet
Laplacian in $\Lsqr(\Omega)$
\begin{equation}
  \label{eq:A0}
  \sfA_0 u \defeq -\Delta u,
  \qquad \dom\sfA_0 = \Sobn^1(\Omega)\cap \Sob^2(\Omega).
\end{equation}
Our second main result concerns the case of non-resonant potentials.
Let us introduce the parameter
\begin{equation}\label{eq:COmega}
	C_\Omega := \inf_{u\in \Sobn^1(\Omega)\setminus\{0\}}\frac{\displaystyle
		\int_\Omega|\nabla u|^2\dd x}{\displaystyle \int_\Omega|u(x)|^2(\dist(x,\Sigma))^{-2}\dd x}.
\end{equation}
By~\cite[Theorem 3.2.1, discussion below Eq. (3.2.4)]{BEL} we have
$C_\Omega = \frac14$ for bounded smooth convex domains,
 while for general bounded $C^\infty$-smooth domains $C_\Omega \in (0,\frac14]$. 
\begin{maintheorem}[the non-resonant case]
  \label{thm:main2}
  Assume that $V\in \Cci(\clo \dR_+)$ is
  such that $\max_{t\in\dR_+} (-t^2 V(t)) < C_\Omega$. Then the family of operators $\sfA_\eps$ converges to
  $\sfA_0$ in strong resolvent sense as $\eps \to 0$.
\end{maintheorem}
It follows from the Hardy inequality on an interval~\cite[Eq. 3.3.4]{BEL} 
that any potential satisfying the assumption in Theorem~\ref{thm:main2} 
is non-resonant, since the Schr\"odinger operator in Remark~\ref{rem:resonant}~\ref{resonant.a} 
is strictly positive and thus can not have eigenvalue zero. Any non-negative potential satisfies 
the assumption in the above theorem. Moreover, a small negative part is allowed. 
The assumption imposes no restriction on the positive part of the potential.
The proof of this theorem relies on the same circle of ideas as the
proof of Theorem~\ref{thm:main}. In the construction of the
identification operators we use the derivative of the non-resonant
solution instead of the non-resonant solution itself as we do in the
proof of Theorem~\ref{thm:main} with the resonant solution.  By
analogy with the one-dimensional case, we expect that $\sfA_\eps$
converge as $\eps \to 0$ to $\sfA_0$ in strong resolvent sense for
general non-resonant potential $V$. However, we have not been able to
find a proof for this claim.

In the proofs of both Theorems~\ref{thm:main} and~\ref{thm:main2}, the
use of identification operators with cut-off functions based on the
solution $\psi_0$ leads to cancellation of `bad' terms; see
Lemmata~\ref{lem:form-diff.a}
and~\ref{lem:form-diff.a'}. Identification operators with the same
mapping properties, but with other cut-off functions, would not lead
to such a cancellation.

\begin{remark}[no uniform ellipticity]
  The condition on  $V$ in the non-resonant case implies that there
  is a constant $c >0$  such that
  $\norm{(\sfA_\eps+\ii)^{-1} v}_{\Sob^1(\Omega)} \le c \| v \|_{\Lsqr(\Omega)}$ for all
  $v \in \Lsqr(\Omega)$ and all $\eps$ small enough (see
  Lemma~\ref{lem:coercive.est}).  If the latter estimate does not
  hold, then we can show that the family of sesquilinear forms associated with the operators $\sfA_\eps$
  is not uniformly elliptic (see Lemma~\ref{lem:no.coercive.est}), a
  concept called ``equi-elliptic'' in~\cite{mnp:13}.
\end{remark}

%
\section{Preliminaries}
\label{sec:prelim}
%

All operators and forms act in the Hilbert space $\Lsqr(\Omega)$; we
denote its norm simply by
$\norm u:=(\int_\Omega \abs {u(x)}^2 \dd x)^{1/2}$.  $\Lsqr$-norms of
subsets $\Omega'\subset \Omega$ and similar norms are typically
indicated by a corresponding subscript such as
$\norm u_{\Lsqr(\Omega')}$.
\subsection{Tubular coordinates}\label{ssec:tubular}
In this subsection, we briefly recall main properties of tubular
coordinates. For any $t > 0$, we will use the notation
$\Omega_t = \{x\in\Omega\colon \dist(x,\Sigma) <
t\}\subset\Omega$, where $\dist(\cdot,\Sigma)$ is the Euclidean distance to $\Sigma$.  By~\cite[Theorem 5.25]{Lee} there exists a
sufficiently small $\delta > 0$ such that the mapping
\begin{equation}
  \label{eq:tub.coord}
  \Phi\colon\Sigma\times (0,\delta)\to\dR^d,
  \qquad \Phi(s,t) \defeq s-t\nu(s)
\end{equation}
is a diffeomorphism onto $\Omega_\delta$. This mapping defines
coordinates $(s,t)$ in $\Omega$ on the tubular neighbourhood
$\Omega_{\delta}$ of $\Sigma$. The metric $G$ induced on
$\Sigma\times (0,\delta)$ by this embedding is
\begin{equation}
  \label{eq:G}
  G = g\circ(\sfI_s - t\sfL_s)^2 + \dd t^2,
\end{equation}
where $\sfI_s\colon T_s\Sigma\to T_s\Sigma$ is the identity map, and
$g$ is the metric on $\Sigma$ induced by the embedding into
$\dR^d$. The volume form associated with the metric $G$ on
$\Sigma\times(0,\delta)$ is given by
\begin{equation}
  \label{eq:vf} %
  \abs{\det G}^{1/2} \dd s \dd t
  = \phi(s,t) \abs{\det g}^{1/2}\dd s\dd t
  = \phi(s,t)\dd\sigma(s)\dd t,
\end{equation}
where
\begin{equation}
  \label{eq:phi}
  \phi\colon\Sigma\times(0,\delta)\rightarrow\dR,\qquad\phi(s,t)
  = \abs*{\det(\sfI_s - t\sfL_s)}
  = 1 - (d-1) H(s)t + p(s,t)t^2
\end{equation}
and where $p$ is a polynomial in $t$ with $\Ci$-smooth coefficients
depending on $s$. We will also make use of the following constant
\begin{equation}
  \label{eq:mu}
  \rho
  \defeq \min_{(s,t)\in \Sigma\times(0,\delta)}\phi(s,t)>0;
\end{equation}
Note that $\rho>0$ is automatically fulfilled as $\Sigma$ is compact,
$\phi$ is continuous, and $\Phi$ is a diffeomorphism.  Let us choose
an orthonormal local coordinate system $(e_1(s),\dots, e_{d-1}(s))$ on
$\Sigma$ at $s\in\Sigma$. By~\eqref{eq:G}, the matrix
$(G_{jk})_{j,k=1}^d$ of the metric $G$ in the local coordinate system
$(e_1(s),\dots,e_{d-1}(s),\nu(s))$ on $\Omega_\delta$ at
$x = \Phi(s,t)$ has block structure and, in particular,
$G_{jd} = G_{dj} = \delta_{jd}$ for all $j\in\{1,2,\dots,d\}$;
\cf~\cite[Lemma 2.3]{LO23}.

Let us define the following unitary map
\begin{equation*}
  \sfU \colon \Lsqr(\Omega_\delta)
  \to \Lsqr(\Sigma\times(0,\delta);\phi(s,t)\dd\sigma(s)\dd t),
  \qquad (\sfU u)(s,t) = u(\Phi(s,t)).
\end{equation*}
For any $u,v \in \Sob^1(\Omega_\delta)$ with the notation
$\wt u \defeq\sfU u$ and $\wt v \defeq \sfU v$, we obtain
\begin{align}
  \nonumber
  \int_{\Omega_\delta} \nabla u\conj{\nabla v}\dd x
  &= \int_0^\delta\int_\Sigma \sum_{j,k=1}^d
    G^{jk}\partial_j \wt u
     \conj{\partial_k \wt v}\phi(s,t)\dd\sigma(s)\dd t\\
  \label{eq:grad}
  &= \int_0^\delta \int_\Sigma
    \bigg(\sum_{j,k=1}^{d-1} G^{jk}\partial_j \wt u \conj{\partial_k \wt v}
    + \partial_d \wt u\conj{\partial_d \wt v}\bigg) \phi(s,t)\dd\sigma(s)\dd t,
\end{align}
where the derivatives $\partial_j$ for $j=1,\dots,d-1$ on $\Sigma$
correspond to the choice of local coordinate system and where $(G^{jk})_{j,k=1}^d$ stands for the inverse to the matrix $(G_{jk})_{j,k=1}^d$.  We also
write $\partial_t$ for $\partial_d$ if we need to stress that the
derivative is with respect to the $d$-th variable $t$.

\subsection{Quadratic forms and operators}
The self-adjoint Robin Laplacian $\sfA$ defined in~\eqref{eq:A} with
mean curvature entering the boundary condition is associated with the
following closed, densely defined, symmetric, and lower-semibounded
quadratic form
\begin{equation*}
  \fra[u] \defeq \norm{\nabla u}^2_{\Lsqr(\Omega;\dC^d)} +
  \frac{d-1}{2}\int_\Sigma H(s)\abs{u(s)}^2\dd\sigma(s),
  \quad \dom\fra \defeq \Sob^1(\Omega)
\end{equation*}
in the Hilbert space $\Lsqr(\Omega)$.  Assume that $\eps > 0$ is so
small such that $\supp V \subset [0,\delta\eps^{-1})$ holds.  Then the
Schr\"odinger operator $\sfA_\eps$ defined in~\eqref{eq:Aeps} is
associated with the closed, densely defined, symmetric, and
lower-semibounded quadratic form
\begin{equation*}
  \fra_\eps[u]
  \defeq \norm{\nabla u}^2_{\Lsqr(\Omega;\dC^d)}
  + \frac{1}{\eps^2}
  \int_0^{\delta}\int_\Sigma V
  \Bigl(\frac{t}{\eps}\Bigr)
  \abs{u(\Phi(s,t))}^2\phi(s,t)\dd \sigma(s)\dd t,
  \quad \dom\fra_\eps \defeq \Sobn^1(\Omega)
\end{equation*}
in the Hilbert space $\Lsqr(\Omega)$.  Finally, the Dirichlet Laplacian
$\sfA_0$ is associated with the closed, non-negative, densely defined
quadratic form in $\Lsqr(\Omega)$ defined by
\begin{equation*}
  \fra_0[u] \defeq \norm{\nabla u}^2_{\Lsqr(\Omega;\dC^d)},
  \quad \dom\fra_0 \defeq \Sobn^1(\Omega).
\end{equation*}
Let us denote the resolvents of $\sfA_\eps$, $\sfA$, and $\sfA_0$ (at
the point $\lambda = \ii$) by
\begin{equation}
  \label{eq:resolvents}
  \sfR_\eps \defeq (\sfA_\eps -\ii)^{-1},
  \qquad
  \sfR \defeq (\sfA-\ii)^{-1},\qquad\sfR_0 \defeq (\sfA_0-\ii)^{-1}.
\end{equation}
Note also that by elliptic regularity~\cite[Theorem 4.18]{McL} for any
$u\in \Cci(\Omega)$ we have
$\sfR_\eps^*u, \sfR u,\sfR_0 u\in \Ci(\ov\Omega)$.  This
observation will be used in the proofs of Theorems~\ref{thm:main}
and~\ref{thm:main2}.

%
\section{Convergence for resonant potentials}
\label{sec:res-pot}
%

We split the proof of Theorem~\ref{thm:main} into several steps.  We
first define some auxiliary boundary mappings which will be used in a
convenient representation for the difference of the sesquilinear forms
of $\sfA$ and $\sfA_\eps$.  In this representation, we also use an
identification operator defined in a second step and mapping functions
from the form domain of $\sfA$ into the form domain of $\sfA_\eps$
(the latter requires a Dirichlet boundary condition on $\bd \Omega$).
The identification operator is basically multiplication with the
scaled resonant solution $\psi_0$ of the
initial-value problem~\eqref{eq:IVP}.

\subsection{Auxiliary boundary mappings}
In this subsection we prove an auxiliary estimate in the neighbourhood
of the boundary.  Let us define for $t\in [0,\delta)$ the mappings
\begin{subequations}
  \label{eq:GUt}
  \begin{align}
    \label{eq:GUt.a}
    \Gamma_t&\colon \dom\sfA_\eps\cap \Ci(\ov\Omega)\to \Lsqr(\Sigma),
    &(\Gamma_t v)(s)
    &\defeq \wt v(s,t),\\
    \label{eq:GUt.b}
    \Upsilon_t &\colon \dom\sfA\cap \Ci(\ov\Omega) \to \Lsqr(\Sigma),
    &(\Upsilon_t u)(s)
    &\defeq 2\phi(s,t)(\partial_t \wt u)(s,t) + \partial_t\phi(s,t)\wt u(s,t),
  \end{align}
\end{subequations}
where the function $\phi$ is as in~\eqref{eq:phi} and where the notation
\begin{equation}
  \label{eq:fcts.coord}
  \wt u
  =\sfU (u\restr{\Omega_\delta})
  = u\circ \Phi
  \qquadtext{and}
  \wt v
  = \sfU (v \restr {\Omega_\delta})
  = v\circ \Phi
\end{equation}
for $u\in\dom\sfA\cap \Ci(\ov\Omega)$ and
$v\in\dom\sfA_\eps\cap \Ci(\ov\Omega)$ is employed.
Note that the functions $\wt u$ and $\wt v$ are
smooth ($\wt u, \wt v \in \Ci(\Sigma\times[0,\delta))$).  Moreover, we
have $\wt v(s,0) =0$ for all $s\in\Sigma$.

The auxiliary mappings $\Gamma_t$ and $\Upsilon_t$ will appear in an
expression for the difference of the limit and approximating
sesquilinear forms, \cf Lemma~\ref{lem:form-diff.a} below.

For an open set $\Omega'\subset\dR^d$ we define the following norms in the
Sobolev spaces $\Sob^1(\Omega')$ and $\Sob^2(\Omega')$
\begin{equation*}
  \norm u ^2_{\Sob^1(\Omega')}
  \defeq \int_{\Omega'}\big(\abs{\nabla u}^2 + \abs u^2\bigr)\dd x,
  \qquad
  \norm u^2_{\Sob^2(\Omega')}
  \defeq \int_{\Omega'}\big(\abs{D^2 u}^2 + \abs{\nabla u}^2 + \abs u^2\bigr)\dd x,
\end{equation*}
where $\abs{D^2u}$ stands for the Hilbert-Schmidt norm of the Hessian
of $u$.

We now estimate one of the auxiliary mappings
\begin{lemma}
  \label{lem:aux}
  Let the mapping $\Upsilon_t$ be defined as in~\eqref{eq:GUt.b} and the
  operator $\sfA$ be as in~\eqref{eq:A}.  Then, there exists a
  constant $c > 0$ such that for any $t\in [0,\delta)$
  \begin{equation*}
    \norm{\Upsilon_t u}_{\Lsqr(\Sigma)}
    \le c \sqrt{t} \norm u_{\Sob^2(\Omega_t)}
  \end{equation*}
  holds for all $u\in \dom\sfA\cap \Ci(\ov\Omega)$.
\end{lemma}
\begin{proof}
  For $u\in\dom\sfA\cap \Ci(\clo \Omega)$ we use again the
  notation~\eqref{eq:fcts.coord}. Combining the boundary
  condition~\eqref{eq:A} together with the identities $\phi(s,0) = 1$
  and $\partial_t\phi(s,0) = - (d-1)H(s)$ for $s \in \Sigma$ we see
  that
  \begin{align*}
    (\Upsilon_0 u)(s)
    &= 2(\partial_t \wt u)(s,0) - (d-1)H(s) \wt u(s,0)\\
    &= -2\bigg(\partial_\nu u\restr{\Sigma}
      + \frac{(d-1)H}{2} u\restr{\Sigma}\bigg)(\Phi(s,0)) = 0.
  \end{align*}
  By the fundamental theorem of calculus and~\eqref{eq:GUt.b}, we
  obtain
  \begin{align*}
    (\Upsilon_t u)(s)
    &=
    \int_0^t \frac{\partial((\Upsilon_tu)(s))}{\partial t}\Big|_{t=t'} \dd t'\\
    &=\int_0^t\bigl(
      2\phi(s,t')\partial_t^2 \wt u(s,t')
      +3\partial_t\phi(s,t')\partial_t \wt u(s,t')
      + \partial_t^2\phi(s,t')\wt u(s,t')
    \bigr)\dd t'
  \end{align*}
  for any $t\in [0,\delta)$ and any $s\in\Sigma$.  In view
  of~\eqref{eq:phi} there exists a constant $C > 0$ such that
  $\abs{\phi(s,t)}$, $\abs{\partial_t\phi(s,t)}$,
  $\abs{\partial_t^2\phi(s,t)} \le C$ for all $s\in\Sigma$ and
  $t\in [0,\delta)$.  Applying the Cauchy-Schwarz inequality we obtain
  \begin{equation}
    \label{eq:bound}
    \abs{(\Upsilon_t u)(s)}
    \le 3\sqrt{3}C\sqrt{t}
    \Bigl(\int_0^t\big(\abs{\partial_t^2 \wt u(s,t')}^2
    + \abs{\partial_t \wt u(s,t')}^2
    + \abs{\wt u(s,t')}^2\bigr) \dd t'\Bigr)^{1/2}.
  \end{equation}
  We have
  \begin{equation*}
    \partial_t \wt u(s,t)
    = -\iprod{\nabla u(\Phi(s,t))}{\nu(s)}_{\dR^d},\qquad
    \partial_t^2 \wt u(s,t)
    = \iprod[\big]{D^2u(\Phi(s,t))\nu(s)}{\nu(s)}_{\dR^d},
  \end{equation*}
  where $\iprod \cdot \cdot_{\dR^d}$ stands for the standard inner
  product in $\dR^d$.  Hence, we obtain
  \begin{align*}
    \norm{\Upsilon_t u}^2_{\Lsqr(\Sigma)}
    &\le 27 C^2t\int_\Sigma\int_0^t
      \big(\abs{D^2 u(\Phi(s,t'))}^2
      + \abs{\nabla u(\Phi(s,t'))}^2 + \abs{u(\Phi(s,t'))}^2
      \bigr)\dd t'\dd \sigma(s)\\
    &\le \frac{27 C^2t}{\rho}
      \int_\Sigma\int_0^t
     \big(\abs{D^2 u(\Phi(s,t'))}^2
     + \abs{\nabla u(\Phi(s,t'))}^2 + \abs{u(\Phi(s,t'))}^2
     \bigr)\phi(s,t') \dd t' \dd \sigma(s)\\
    &\le \frac{27 C^2t}{\rho}\norm u^2_{\Sob^2(\Omega_t)},
  \end{align*}
  where the constant $\rho$ is as in~\eqref{eq:mu}.  Hence, the
  inequality in the formulation of the lemma holds with
  $c = (3\sqrt{3}C)/\sqrt \rho$.
\end{proof}

\subsection{The identification   operator and an expression for the form
  difference}

For $\eps > 0$, we define the self-adjoint bounded multiplication
operator
\begin{equation*}
  \sfJ_\eps \colon \Lsqr(\Omega)\to \Lsqr(\Omega),
  \qquad
  (\sfJ_\eps u)(x)
  \defeq \psi_0\Bigl(\frac{\dist(x,\Sigma)}{\eps}\Bigr) u(x),
\end{equation*}
where $\psi_0\in \Ci(\clo \dR_+)$ is the bounded solution of the
initial-value problem~\eqref{eq:IVP} satisfying $\psi_0(t) =1$ for all
$t > a$.

It follows from
$\Ci$-smoothness of the mapping
$\dist(\cdot\,,\Sigma) \colon \Omega_\delta \to \dR_+$ and of the function
$\psi_0$ that for all sufficiently small $\eps > 0$ it holds that
$\ran(\sfJ_\eps\restr{\Sob^1(\Omega)})\subset \Sobn^1(\Omega)$ and
$\ran(\sfJ_\eps\restr{\Ci(\clo\Omega)})\subset
\Ci(\clo\Omega)$, where for the first-mentioned property we took
into account $\psi_0(0) = 0$.

We first compare the identification operator $\sfJ_\eps$ with the
identity $\sfI$:
\begin{lemma}
  \label{lem:j-id}
  For any $u \in \Lsqr(\Omega)$ we have
  \begin{equation*}
    \norm{(\sfJ_\eps  -\sfI)u}
    \le \norm{\psi_0-1}_\infty \norm u_{\Lsqr(\Omega_{a\eps})}
    \to 0 \qquadtext{as}\eps\to 0.
  \end{equation*}
  In other words, $\sfJ_\eps$ converges to the identity operator in
  strong operator sense as $\eps\to0$.
\end{lemma}
\begin{proof}
  We actually have
  \begin{equation*}
    \norm{(\sfJ_\eps  -\sfI)u}^2
    = \int_\Omega
    \abs*{\psi_0\Bigl(\frac{\dist(x,\Sigma)}{\eps}\Bigr)-1}^2 \abs{u(x)}^2\dd x
  \end{equation*}
  from which the desired inequality follows.
\end{proof}

We now see the reason for defining the auxiliary boundary mappings and
the choice of identification operator:
\begin{lemma}
  \label{lem:form-diff.a}
  Assume that $\eps < a^{-1}\delta$.
  Then we have
  \begin{equation*}
    \fra[u,\sfJ_\eps v]
    - \fra_\eps[\sfJ_\eps u, v]
    = \frac{1}{\eps}\int_0^{a\eps}\psi_0'\Bigl(\frac{t}{\eps}\Bigr)
    \iprod[\big]{\Upsilon_t u}{\Gamma_t v}_{\Lsqr(\Sigma)} \dd t
  \end{equation*}
  for $u \in \dom \sfA\cap \Ci(\clo\Omega)$ and
  $v \in\dom \sfA_\eps\cap \Ci(\clo\Omega)$.
\end{lemma}
\begin{proof}
Recall that $\delta > 0$ is chosen as in Subsection~\ref{ssec:tubular} so that the mapping $\Phi$ in \eqref{eq:tub.coord} is a diffeomorphism from $\Sigma\times(0,\delta)$ onto $\delta$-neighbourhood $\Omega_\delta$ of the boundary $\partial\Omega$.
  Under the assumption $\eps < a^{-1}\delta$ the tubular
  coordinates~\eqref{eq:tub.coord} can be used.  In particular, using
  the notation~\eqref{eq:fcts.coord}, we have
  $\wt u, \wt v \in \Lsqr(\Sigma\times
  (0,\delta);\phi(s,t)\dd\sigma(s)\dd t)$ and these functions are $C^\infty$-smooth.

  Clearly, the contributions of $\fra[u,\sfJ_\eps v]$ and
  $\fra_\eps[\sfJ_\eps u, v]$ outside the tubular neighbourhood
  $\Omega_{a\eps}$ cancel.  Using Equation~\eqref{eq:grad} we also
  observe that the contribution of the gradient terms corresponding to
  derivatives in the direction tangential to $\Sigma$ cancel too.  We
  end up with the following formula
  \begin{align*}
    \fra[u,\sfJ_\eps v]
    -
    \fra_\eps[\sfJ_\eps u, v]
    &=
      \int_\Sigma\int_0^{a\eps}
      \partial_t\wt u(s,t)
      \partial_t\Bigl(
      \psi_0\Bigl(\frac{t}{\eps}\Bigr) \conj{\wt v(s,t)}\Bigr)
      \phi(s,t)\dd t\dd \sigma(s)\\
    &\qquad\qquad
      -
      \int_\Sigma\int_0^{a\eps}
      \partial_t\Bigl(
      \psi_0\Bigl(\frac{t}{\eps}\Bigr) \wt u(s,t)\Bigr)
      \conj{\partial_t \wt v(s,t)}
      \phi(s,t)\dd t\dd \sigma(s)\\
    &\qquad\qquad-\frac{1}{\eps^2}\int_\Sigma\int_0^{a\eps}
      V\Bigl(\frac{t}{\eps}\Bigr)\psi_0\Bigl(\frac{t}{\eps}\Bigr)
      \wt u(s,t) \conj{\wt v(s,t)}\phi(s,t)\dd t\dd \sigma(s) \\
    &=
      \frac{1}{\eps}\int_\Sigma\int_0^{a\eps} \partial_t \wt u(s,t)
      \psi_0'\Bigl(\frac{t}{\eps}\Bigr) \conj{\wt v(s,t)}\phi(s,t)\dd t\dd \sigma(s)\\
    &\qquad\qquad -
      \frac{1}{\eps}\int_\Sigma\int_0^{a\eps}
      \psi_0'\Bigl(\frac{t}{\eps}\Bigr)
      \wt u(s,t)\conj{\partial_t \wt v(s,t)}\phi(s,t)\dd t\dd \sigma(s)\\
    &\qquad\qquad -\frac{1}{\eps^2}\int_\Sigma\int_0^{a\eps}
      V\Bigl(\frac{t}{\eps}\Bigr)\psi_0
      \Bigl(\frac{t}{\eps}\Bigr) \wt u(s,t) \conj{\wt v(s,t)}\phi(s,t)\dd t\dd \sigma(s),
  \end{align*}
  where in the second step two terms cancelled upon using the product
  rule for differentiation.  After integration by parts in the second
  term on the right hand side of the above formula, we arrive at
  \begin{align}
    \nonumber
    \fra[u,\sfJ_\eps v]
    - \fra_\eps[\sfJ_\eps u, v]
    &= \frac{1}{\eps^2}\int_\Sigma\int_0^{a\eps}
      \psi_0''\Bigl(\frac{t}{\eps}\Bigr)
      \wt u(s,t)\conj{\wt v(s,t)}\phi(s,t)\dd t\dd\sigma(s)\\
    \nonumber
    &\qquad +\frac{1}{\eps}\int_\Sigma\int_0^{a\eps}
      \psi_0'\Bigl(\frac{t}{\eps}\Bigr)\wt u(s,t)\conj{\wt v(s,t)}
      \partial_t\phi(s,t)\dd t\dd\sigma(s)\\
    \nonumber
    &\qquad +\frac{2}{\eps} \int_\Sigma\int_0^{a\eps}\psi_0'
      \Bigl(\frac{t}{\eps}\Bigr)
      \partial_t \wt u(s,t)\conj{\wt v(s,t)}\phi(s,t)\dd t\dd \sigma(s)\\
    \label{eq:form-diff.a1}
    &\qquad -\frac{1}{\eps^2}\int_\Sigma\int_0^{a\eps}
      V\Bigl(\frac{t}{\eps}\Bigr)\psi_0\Bigl(\frac{t}{\eps}\Bigr) \wt u(s,t)
      \conj{\wt v(s,t)}\phi(s,t)\dd t\dd\sigma(s),
  \end{align}
  where the boundary term at $t = 0$ vanishes due to $\wt v(s,0) = 0$
  while the boundary term at $t = a\eps$ vanishes due to
  $\psi_0'(a) = 0$.  Using that $\psi_0$ satisfies the differential
  equation $-\psi_0'' + V\psi_0 = 0$ we note that the first and the
  last terms on the right hand side in the above formula cancel each
  other.  The remaining two integrals just give the desired expression
  involving $\Gamma_t$ and $\Upsilon_t$.
\end{proof}

We now estimate the expression of Lemma~\ref{lem:form-diff.a}:
\begin{lemma}
  \label{lem:form-diff.a.est}
  Assume that $\eps < a^{-1}\delta$.
  Then we have
  \begin{equation*}
    \abs[\big]{\fra[u,\sfJ_\eps v]
      - \fra_\eps[\sfJ_\eps u, v]}
    \le \hat c \norm u_{\Sob^2(\Omega_{a\eps})}
  \norm v_{\Lsqr(\Omega_{a\eps})}
  \end{equation*}
  for $u \in \dom \sfA\cap \Ci(\clo\Omega)$ and
  $v \in \dom \sfA_\eps\cap \Ci(\clo\Omega)$, where $\hat c$ is given
  in~\eqref{eq:def.c.hat} below.
\end{lemma}
\begin{proof}
  Using Cauchy-Schwarz inequality (twice), Lemma~\ref{lem:aux} and
  Lemma~\ref{lem:form-diff.a} we obtain
  \begin{align}
    \nonumber
    \abs*{\fra[u,\sfJ_\eps v] - \fra_\eps[\sfJ_\eps u, v]}
    &\le \frac{\norm{\psi_0'}_{\infty}}{\eps}
      \int_0^{a\eps}
      \norm{\Upsilon_t u}_{\Lsqr(\Sigma)}
      \norm{\Gamma_t v}_{\Lsqr(\Sigma)}\dd t\\
    \nonumber
    &\le
      \frac{c \norm{\psi_0'}_{\infty}\sqrt a}{\sqrt\eps}
      \norm{u}_{\Sob^2(\Omega_{a\eps})}
      \int_0^{a\eps}
      \norm{\Gamma_t v}_{\Lsqr(\Sigma)}\dd t\\
    \nonumber
    &\le
      ca\norm{\psi_0'}_{\infty}
      \norm{u}_{\Sob^2(\Omega_{a\eps})}
      \Bigl(\int_0^{a\eps}\int_\Sigma
      \abs{\wt v(s,t)}^2 \dd\sigma(s)\dd t\Bigr)^{1/2}\\
    \label{eq:def.c.hat}
    &\le \hat c
      \norm{u}_{\Sob^2(\Omega_{a\eps})}
      \norm{v}_{\Lsqr(\Omega_{a\eps})}
      \qquadtext{where}
      \hat c \defeq \frac{ca\norm{\psi_0'}_{\infty}}{\sqrt\rho},
  \end{align}
  using also~\eqref{eq:mu} for the last estimate.
\end{proof}

\subsection{Proof of Theorem~\ref{thm:main}}

\begin{proof}[Proof of Theorem~\ref{thm:main}]
  For any $u,v\in \Cci(\Omega)$ we obtain
  \begin{align*}
    \iprod[\big]{(\sfR_\eps-\sfR)u} v
    &= \iprod[\big] u {\sfJ_\eps\sfR_\eps^*v}
    - \iprod[\big]{\sfJ_\eps\sfR u} v\\
    &\qquad\qquad+ \iprod[\big] u{ (\sfI-\sfJ_\eps)\sfR_\eps^* v}
      -
     \iprod[\big] {(\sfI-\sfJ_\eps)\sfR u} v\\
    &= \iprod[\big] {(\sfA-\ii)\sfR u}{\sfJ_\eps\sfR_\eps^* v}
      - \iprod[\big]{\sfJ_\eps\sfR u}
         {(\sfA_\eps +\ii)\sfR_\eps^* v}\\
    &\qquad\qquad
      + \iprod[\big] u {(\sfI-\sfJ_\eps)\sfR_\eps^* v}
      - \iprod[\big]{(\sfI-\sfJ_\eps)\sfR u} v\\
    &=
      \fra[\sfR u, \sfJ_\eps\sfR_\eps^* v]
      - \fra_\eps[\sfJ_\eps\sfR u,\sfR ^*_\eps v]\\
    &\qquad\qquad
      + \iprod[\big] {(\sfI-\sfJ_\eps)u} {\sfR_\eps^* v}
      - \iprod[\big]{(\sfI-\sfJ_\eps)\sfR u} v,
  \end{align*}
  using first representation theorem for the sesquilinear forms $\fra$
  and $\fra_\eps$ associated with $\sfA$ and $\sfA_\eps$,
  respectively, and using also the self-adjointness of
  $\sfI-\sfJ_\eps$ for the last equality.  Moreover, we have
  $\norm{\sfR u} \le \norm u$ and $\norm{\sfR_\eps^*v} \le \norm v$ as
  $\sfR$ and $\sfR_\eps^*$ are the resolvents at the points $\pm\ii$
  and $\sfA$ and $\sfA_\eps$ are self-adjoint, respectively, hence we
  conclude (using Cauchy-Schwarz)
  \begin{align*}
    \abs[\big]{\iprod[\big]{(\sfR_\eps-\sfR)u} v%
    }
    &\le \abs[\big]{\fra[\sfR u, \sfJ_\eps\sfR_\eps^* v]
      - \fra_\eps[\sfJ_\eps\sfR u,\sfR ^*_\eps v]}
    + \norm[\big] {(\sfJ_\eps-\sfI)u}\norm{\sfR_\eps^*v}
    + \norm[\big] {(\sfJ_\eps-\sfI)\sfR u}
    \norm v\\
    &\le \Bigl(
        \hat c \norm{\sfR u}_{\Sob^2(\Omega_{a\eps})}
        + \norm{\psi_0-1}_\infty
      \bigl(
        \norm u_{\Lsqr(\Omega_{a\eps})}
        +\norm {\sfR u}_{\Lsqr(\Omega_{a\eps})}
      \bigr)
      \Bigr)\norm v
  \end{align*}
  using also Lemmas~\ref{lem:j-id} and~\ref{lem:form-diff.a.est} for
  the second estimate (note that
  $\sfR u,\sfR_\eps^* v\in \Ci(\clo\Omega)$).  From the
  characterisation of the dual of a Hilbert space on the dense subset
  $\Cci(\Omega)$ of $\Lsqr(\Omega)$, we obtain
  \begin{align*}
    \norm{\sfR_\eps u -\sfR u}
    &=\sup_{\substack{v\in \Cci(\Omega)\\ \norm v = 1}}
    \abs[\Big]{\iprod[\big]{(\sfR_\eps - \sfR)u} v}\\
    & \le
        \hat c \norm{\sfR u}_{\Sob^2(\Omega_{a\eps})}
        + \norm{\psi_0-1}_\infty
      \bigl(
        \norm u_{\Lsqr(\Omega_{a\eps})}
        +\norm {\sfR u}_{\Lsqr(\Omega_{a\eps})}
      \bigr).
  \end{align*}
  Now all norms on $\Omega_{a\eps}$ converge to $0$ by Lebesgue's
  convergence theorem for $u \in \Cci(\Omega)$.  By density of
  $\Cci(\Omega)$ in $\Lsqr(\Omega)$, we conclude that $\sfA_\eps$
  converges to $\sfA$ in strong resolvent sense.
\end{proof}

%
\section{Convergence for non-resonant potentials}
\label{sec:non-neg.pot}

\subsection{Auxiliary boundary mappings, the identification operator
  and some related estimates}

In this subsection, we provide another lemma needed in the proof of
Theorem~\ref{thm:main2}.  Recall that the mapping $\Gamma_t$ is
defined in~\eqref{eq:GUt}.
%
\begin{lemma}
  \label{lem:aux2}
  We have
  \begin{equation*}
    \norm{\Gamma_tu}_{\Lsqr(\Sigma)}
    \le \sqrt{\frac t\rho}\norm{\nabla u}_{\Lsqr(\Omega_t;\dC^d)}
  \end{equation*}
  for all $u\in \dom\sfA_\eps\cap C^\infty(\ov\Omega)$ and $t\in(0,\delta)$.
\end{lemma}
\begin{proof}
  Let $t\in(0,\delta)$ and $u\in \dom\sfA_\eps\cap C^\infty(\ov\Omega)$.  By the fundamental
  theorem of calculus we obtain in view of $\Gamma_0 u = 0$ that
  \begin{equation*}
    (\Gamma_t u)(s)
    = \int_0^t\partial_t \wt u(s,t')\dd t'
  \end{equation*}
  for any $t\in[0,\delta)$ and $s\in\Sigma$ (recall that
  $\wt u=u \circ \Phi$ as in~\eqref{eq:fcts.coord}).  Using the
  Cauchy-Schwarz inequality we obtain
  \begin{equation*}
    \abs{(\Gamma_t u)(s)}
    \le \sqrt{t}
    \Bigl(\int_0^t\abs{\partial_t \wt u(s,t')}^2\dd t'\Bigr)^{1/2}.
  \end{equation*}
  As
  $\partial_t \wt u(s,t) = -\iprod {\nabla
    u(\Phi(s,t))}{\nu(s)}_{\dR^d}$, we deduce (using Cauchy-Schwarz
  again)
  \begin{align*}
    \norm{\Gamma_t u}^2_{\Lsqr(\Sigma)}
    &\le t\int_\Sigma\int_0^t\abs{\nabla u(\Phi(s,t'))}^2
      \dd t'\dd\sigma(s)\\
    & \le\frac{t}{\rho}
      \int_\Sigma\int_0^t \abs{\nabla u(\Phi(s,t'))}^2\phi(s,t')
      \dd t' \dd\sigma(s)\\
    & = \frac{t}{\rho}
      \int_{\Omega_t}\abs{\nabla u}^2\dd x
      =\frac{t}{\rho} \norm{\nabla u}^2_{\Lsqr(\Omega_t;\dC^d)}.
    \qedhere
  \end{align*}
\end{proof}

We also need the following modified boundary mapping
$\wt \Gamma_t \colon \dom\sfA_\eps\cap C^\infty(\ov\Omega) \to
\Lsqr(\Sigma)$ defined for all $t\in(0,\delta)$ by
\begin{equation}
  \label{eq:bd.map'}
  (\wt \Gamma_t v)(s)
  \defeq \phi(s,t) \wt v(s,t)
  =\phi(s,t) v(\Phi(s,t)).
\end{equation}
\begin{corollary}
  \label{cor:aux2}
  We have
  \begin{equation*}
    \norm{\wt \Gamma_tv}_{\Lsqr(\Sigma)}
    \le \norm \phi_\infty\sqrt{\frac {t}\rho}\norm{\nabla v}_{\Lsqr(\Omega_t;\dC^d)}
  \end{equation*}
  for all $v\in \dom\sfA_\eps\cap C^\infty(\ov\Omega)$ and $t\in(0,\delta)$.
\end{corollary}
\begin{proof}
  We just estimate $0<\phi(s,t) \le \norm \phi_\infty$ in the first
  step and then use Lemma~\ref{lem:aux2}.
\end{proof}

Recall that we fix $a > 0$ such that $\supp V \subset [0,a]$.  Let
$\psi_0$ be a non-resonant solution of~\eqref{eq:IVP} normalised so
that $\psi_0'(t) = 1$ for all $t > a$; see
Remark~\ref{rem:resonant}~\itemref{resonant.b}.

For $\eps > 0$ we define the self-adjoint bounded multiplication
operator
\begin{equation*}
  \sfK_\eps\colon \Lsqr(\Omega)\to \Lsqr(\Omega),\qquad
  (\sfK_\eps u)(x)
  =
  \psi_0'\Bigl(\frac{\dist(x,\Sigma)}{\eps}\Bigr) u(x).
\end{equation*}
The only difference with $\sfJ_\eps$ is that we use $\psi_0'$ (with
$\psi_0'=1$ outside $(0,a)$) instead of $\psi_0$. We obtain as in the
proof of Lemma~\ref{lem:j-id}:
\begin{lemma}
  \label{lem:j-id'}
  For any $u \in \Lsqr(\Omega)$ we have
  \begin{equation*}
    \norm{(\sfK_\eps  -\sfI)u}
    \le \norm{\psi_0'-1}_\infty \norm u_{\Lsqr(\Omega_{a\eps})}
    \to 0 \qquadtext{as}\eps\to 0.
  \end{equation*}
  In other words, $\sfK_\eps$ converges to the identity operator in
  strong operator sense as $\eps\to0$.
\end{lemma}

For the form difference, we have a similar expression as in
Lemma~\ref{lem:form-diff.a}:
\begin{lemma}
  \label{lem:form-diff.a'}
  Assume that $\eps < a^{-1}\delta$.
  Then we have
  \begin{multline*}
    \fra_0[u,\sfK_\eps v]
    - \fra_\eps[\sfK_\eps u, v]
    = \frac{1}{\eps}\int_0^{a\eps}\psi_0''\Bigl(\frac t \eps\Bigr)
    \iprod[\big]{\Upsilon_t u}{\Gamma_t v}_{\Lsqr(\Sigma)} \dd t\\
    + \frac1{\eps^2}\int_0^{a\eps}
    (V' \psi_0)\Bigl(\frac t \eps\Bigr)
    \iprod[\big]{\Gamma_t u}{\wt \Gamma_t v}_{\Lsqr(\Sigma)} \dd t
  \end{multline*}
  for $u \in \dom \sfA_0\cap C^\infty(\ov\Omega)$ and $v \in \dom \sfA_\eps\cap C^\infty(\ov\Omega)$.
\end{lemma}
\begin{proof}
  The proof is exactly the same as the proof of
  Lemma~\ref{lem:form-diff.a} --- except that by differentiating the
  identity $-\psi_0'' + V\psi_0 = 0$ we obtain
  $-\psi_0''' + V'\psi_0 + V\psi_0' = 0$.  In particular, the first
  and fourth term in~\eqref{eq:form-diff.a1} do not cancel, as now
  $\psi_0'''-V\psi_0'= V'\psi_0$ remains.  In particular, we have
  \begin{align*}
    \fra_0[u,\sfK_\eps v]
    -
    \fra_\eps[\sfK_\eps u, v]
    &=\frac{1}{\eps}\int_\Sigma\int_0^{a\eps}
      \psi_0''\Bigl(\frac{t}{\eps}\Bigr)\wt u(s,t)\conj{\wt v(s,t)}
      \partial_t\phi(s,t)\dd t\dd\sigma(s)\\
    &\qquad +\frac{2}{\eps} \int_\Sigma\int_0^{a\eps}\psi_0''
      \Bigl(\frac{t}{\eps}\Bigr)
      \partial_t \wt u(s,t)\conj{\wt v(s,t)}\phi(s,t)\dd t\dd \sigma(s)\\
    &\qquad +\frac{1}{\eps^2}\int_\Sigma\int_0^{a\eps}
      (\psi_0'''-V \psi_0')\Bigl(\frac{t}{\eps}\Bigr) \wt u(s,t)
      \conj{\wt v(s,t)}\phi(s,t)\dd t\dd\sigma(s),
  \end{align*}
  from which the desired formula follows.  Here the boundary term at
  $t = 0$ vanishes due to $\wt v(s,0) = 0$ while the boundary term at
  $t = a\eps$ vanishes due to $\psi_0''(a) = 0$.
\end{proof}

As before, we now estimate the expression of
Lemma~\ref{lem:form-diff.a'}:
\begin{lemma}
  \label{lem:form-diff.a'.est}
  Assume that $\eps < a^{-1}\delta$. Then we have
  \begin{equation*}
    \abs[\big]{\fra_0[u,\sfK_\eps v]
      - \fra_\eps[\sfK_\eps u, v]}
    \le \hat c_0 \norm u_{\Sob^1(\Omega_{a\eps})}
    \norm v_{\Sob^1(\Omega_{a\eps})}
  \end{equation*}
  for $u \in \dom \sfA_0\cap C^\infty(\ov\Omega)$ and
  $v \in \dom \sfA_\eps\cap C^\infty(\ov\Omega)$, where $\hat c_0$ is
  given in~\eqref{eq:def.c0.hat} below.
\end{lemma}
\begin{proof}
  We estimate the first term in Lemma~\ref{lem:form-diff.a'} as
  follows
  \begin{align}
  	\nonumber
  \bigg|	\frac{1}{\eps}&\int_0^{a\eps}\psi_0''\Bigl(\frac t \eps\Bigr)
  	\iprod[\big]{\Upsilon_t u}{\Gamma_t v}_{\Lsqr(\Sigma)} \dd t\bigg|\\
  	\nonumber
  	&\! \le \frac{\|\psi_0''\|_\infty}{\eps}\left(\int_0^{a\eps}\|\Upsilon_t u\|^2_{\Lsqr(\Sigma)}\dd t\right)^{1/2}\left(\int_0^{a\eps}\|\Gamma_t v\|^2_{\Lsqr(\Sigma)}\dd t\right)^{1/2}\\
  	&\!
  	\nonumber\le
  	\frac{a\|\psi_0''\|_\infty\|v\|_{\Sob^1(\Omega_{a\eps})}}{\sqrt{2\rho}}\left(
  	\int_0^{a\eps}\int_\Sigma\left(8\varphi^2(s,t)|\partial_t \wt u(s,t)|^2 + 2(\partial_t\varphi(s,t))^2|\wt u(s,t)|^2\right)\dd \sigma(s)\dd t\right)^{1/2}
  	\\
  	  	&\!
  	\nonumber\le
  	\frac{2a\|\psi_0''\|_\infty\|v\|_{\Sob^1(\Omega_{a\eps})}}{\rho}\max\{\|\varphi\|_\infty, \|\partial_t\varphi\|_\infty\}\!
  		\left(
  	\int_0^{a\eps}\!\int_\Sigma\left[|\partial_t \wt u(s,t)|^2 + |\wt u(s,t)|^2\right]\varphi(s,t)\dd \sigma(s)\dd t\right)^{1/2}\\
  	&\!
  	\nonumber\le
  	\hat{c}'
  	\|u\|_{\Sob^1(\Omega_{a\eps})}
  	\|v\|_{\Sob^1(\Omega_{a\eps})},\qquad\quad\text{where}\qquad \hat{c}' := 
  	\frac{2a\|\psi_0''\|_\infty}{\rho}
  	\max\big\{\|\varphi\|_\infty, \|\partial_t\varphi\|_\infty\big\},
   \end{align}
	where we used Cauchy-Schwarz inequality in the first estimate and Lemma~\ref{lem:aux2} in the second.
	Employing the above estimate we get		 
  \begin{align}
    \nonumber
    \abs[\big]{\fra_0[&u,\sfK_\eps v] - \fra_\eps[\sfK_\eps u, v]}\\
    \nonumber
    &\le \hat c'
      \norm{u}_{\Sob^1(\Omega_{a\eps})}
      \norm{v}_{\Sob^1(\Omega_{a\eps})}
      + \frac{\norm {V'}_\infty \norm {\psi_0}_{\Linfty((0,a))}}{\eps^2}
      \int_0^{a\eps}
      \norm{\Gamma_t u}_{\Lsqr(\Sigma)}
      \norm{\wt \Gamma_t v}_{\Lsqr(\Sigma)}\dd t\\
    \nonumber
    &\le \hat c'
      \norm{u}_{\Sob^1(\Omega_{a\eps})}
      \norm{v}_{\Sob^1(\Omega_{a\eps})}
      + \frac{\norm {V'}_\infty \norm {\psi_0}_{\Linfty((0,a))}
       \norm \phi_\infty}
        {\rho \eps^2}
      \int_0^{a\eps} \!\! t\dd t
      \norm{\nabla u}_{\Lsqr(\Omega_{a\eps};\dC^d)}
      \norm{\nabla v}_{\Lsqr(\Omega_{a\eps};\dC^d)}\\
    \nonumber
    &= \hat c'
      \norm{u}_{\Sob^1(\Omega_{a\eps})}
      \norm{v}_{\Sob^1(\Omega_{a\eps})}
      + \frac{a^2\norm {V'}_\infty \norm {\psi_0}_{\Linfty((0,a))}
      \norm \phi_\infty}
      {2\rho}
      \norm{\nabla u}_{\Lsqr(\Omega_{a\eps};\dC^d)}
      \norm{\nabla v}_{\Lsqr(\Omega_{a\eps};\dC^d)}\\
    \nonumber
    &\le \Bigl(
      \hat c' + \frac{a^2\norm {V'}_\infty \norm {\psi_0}_{\Linfty((0,a))}
       \norm \phi_\infty}
      {2\rho}
      \Bigr)
      \norm{u}_{\Sob^1(\Omega_{a\eps})}
      \norm{v}_{\Sob^1(\Omega_{a\eps})}\\
    \label{eq:def.c0.hat}
    &\le \hat c_0
      \norm{u}_{\Sob^1(\Omega_{a\eps})}
      \norm{v}_{\Sob^1(\Omega_{a\eps})}
      \qquadtext{where}
      \hat c_0 \defeq
      \Bigl(
      \hat c' + \frac{a^2\norm {V'}_\infty \norm {\psi_0}_{\Linfty((0,a))}
        \norm \phi_\infty}
      {2\rho}
      \Bigr)
  \end{align}
  using Lemma~\ref{lem:aux2} and Corollary~\ref{cor:aux2} in the second estimate.
\end{proof}

\subsection{Proof of Theorem~\ref{thm:main2}}

Before providing the proof of Theorem~\ref{thm:main2}, we need one
more estimate: here, it is essential that the potential 
satisfies an assumption on the negative part:

\begin{lemma}
  \label{lem:coercive.est}
Let the constant $C_\Omega > 0$ be defined as in~\eqref{eq:COmega}. Under the assumption 
  \begin{equation}\label{eq:conditionV}
  C_V := \max_{t\in\dR_+} (-V(t)t^2) < C_\Omega
  \end{equation}
  we have
  \begin{equation*}
    \norm {\sfR_\eps^* v}_{\Sob^1(\Omega)}
    \le \sqrt{c_{\Omega,V}+1}\norm v
  \end{equation*}
  for all $v \in \Lsqr(\Omega)$ and the constant $c_{\Omega,V} 
 := \frac12\big(1-\frac{C_V}{C_\Omega}\big)^{-1} \ge \frac12$. 
\end{lemma}
\begin{proof}
  The estimate $\norm {\sfR_\eps^* v}_{\Lsqr(\Omega)} \le \norm v$
  is clear by spectral calculus.
 Under the assumption we stated on the potential $V$,
  we have the inequality
  \[
  		V_\eps(x) \ge -\frac{C_V}{({\rm dist}\,(x,\Sigma))^2},
  \]
  for all $\eps > 0$ and all $x\in\Omega$.
	We have by the definition of the constant $C_\Omega$ that for all $\hat v\in\dom\sfA_\eps$
  \begin{align*}
    \norm{\nabla \hat v}_{\Lsqr(\Omega;\dC^d)}^2
    &\le 2c_{\Omega,V}\left( \norm{\nabla \hat v}_{\Lsqr(\Omega;\dC^d)}^2
    + \iprod[\Big]{V_\eps \hat v}{\hat v}\right)
      =2c_{\Omega,V}\fra_\eps[\hat v]\\
    &=2c_{\Omega,V}
    	\iprod{\sfA_\eps \hat v}{\hat v}
    \le 
    c_{\Omega,V}\left(\norm{\sfA_\eps \hat v}^2 + \norm{\hat v}^2\right)
    =c_{\Omega,V}\norm{(\sfA_\eps+\ii)\hat v}^2
  \end{align*}
  and using the Cauchy-Schwarz inequality combined with standard inequality $2ab\le a^2+b^2$ for $a,b\ge0$ in the penultimate step.  The desired estimate follows by setting
  $\hat v=\sfR_\eps^* v$.
\end{proof}

\begin{proof}[Proof of Theorem~\ref{thm:main2}]
  For any $u,v\in C^\infty_c(\Omega)$ we obtain as in the proof of
  Theorem~\ref{thm:main} the decomposition
  \begin{align*}
    \iprod[\big]{(\sfR_\eps-\sfR_0)u} v
    &= \fra_0[\sfR_0 u, \sfK_\eps\sfR_\eps^* v]
    - \fra_\eps[\sfK_\eps\sfR_0 u,\sfR ^*_\eps v]
    + \iprod[\big] u {(\sfI-\sfK_\eps)\sfR_\eps^* v}
    - \iprod[\big] {(\sfI-\sfK_\eps)\sfR_0 u} v\\
    &= \fra_0[\sfR_0 u, \sfK_\eps\sfR_\eps^* v]
    - \fra_\eps[\sfK_\eps\sfR_0 u,\sfR ^*_\eps v]
    + \iprod[\big] {(\sfI-\sfK_\eps) u} {\sfR_\eps^* v}
    - \iprod[\big] {(\sfI-\sfK_\eps)\sfR_0 u} v
  \end{align*}
  where we used again that $\sfI-\sfK_\eps$ is self-adjoint.  We now
  have (using Cauchy-Schwarz)
  \begin{align*}
    \abs[\big]{\iprod{(\sfR_\eps - \sfR_0) u} v}
    \hspace*{-10ex}&\\
    &\le \abs[\big]{\fra_0[\sfR_0 u,\sfK_\eps\sfR_\eps^* v]
        - \fra_\eps[\sfK_\eps\sfR_0 u,\sfR^*_\eps v]}
    + \norm{(\sfI-\sfK_\eps)u} \norm {\sfR_\eps^* v}
      + \norm{(\sfI-\sfK_\eps)\sfR_0 u} \norm v\\
    &\le  \hat c_0 \norm {\sfR_0 u}_{\Sob^1(\Omega_{a\eps})}
      \norm {\sfR_\eps^* v}_{\Sob^1(\Omega_{a\eps})}
      + \norm{\psi_0'-1}_\infty
      \bigl(
        \norm u_{\Lsqr(\Omega_{a\eps})}\norm {\sfR_\eps^* v}
        +\norm {\sfR_0 u}_{\Lsqr(\Omega_{a\eps})}\norm v
      \bigr)\\
    &\le \Bigl(\sqrt{c_{\Omega,V}+1}\hat c_0 \norm {\sfR_0 u}_{\Sob^1(\Omega_{a\eps})}
      + \norm{\psi_0'-1}_\infty
      \bigl(
        \norm u_{\Lsqr(\Omega_{a\eps})}
        +\norm {\sfR_0 u}_{\Lsqr(\Omega_{a\eps})}
      \bigr)\Bigr) \norm v
  \end{align*}
  using Lemma~\ref{lem:j-id'} and Lemma~\ref{lem:form-diff.a'.est} in
  the second estimate, and
  $\norm{\sfR_\eps^*v}_{H^1(\Omega)} \le\sqrt{c_{\Omega,V}+1} \norm v$ resp.\
  Lemma~\ref{lem:coercive.est} in the last one.  We conclude that
  \begin{align*}
    \norm{\sfR_\eps u -\sfR_0 u}
    &= \sup_{\substack{v\in \Cci(\Omega)\\ \norm v = 1}}
    \abs[\big]{\iprod{(\sfR_\eps - \sfR_0) u} v}\\
    & \le \sqrt{c_{\Omega,V}+1}\hat c_0 \norm {\sfR_0 u}_{\Sob^1(\Omega_{a\eps})}
      + \norm{\psi_0'-1}_\infty
      \bigl(
      \norm u_{\Lsqr(\Omega_{a\eps})}
      +\norm {\sfR_0 u}_{\Lsqr(\Omega_{a\eps})}
      \bigr).
  \end{align*}
  As before, all norms on $\Omega_{a\eps}$ converge to $0$ by
  Lebesgue's convergence theorem for $u \in \Cci(\Omega)$.  as
  $\eps \to 0$.  By density of $\Cci(\Omega)$ in the respective
  spaces, we conclude that $\sfA_\eps$ converges to $\sfA_0$ in the
  strong resolvent sense.
\end{proof}

\subsection{No uniform ellipticity}
\label{ssec:no-uni-ell}
We finally show that if the estimate in Lemma~\ref{lem:coercive.est} with any constant is not true then
$(\fra_\eps)_\eps$ is not uniformly elliptic:
\begin{lemma}
  \label{lem:no.coercive.est}
  If the estimate in Lemma~\ref{lem:coercive.est} does not hold for
  any constant then
  $(\fra_\eps)_\eps$ is not uniformly elliptic, i.e.,
  \begin{equation*}
    \exists \alpha>0\;
    \exists \omega \in \R\;
    \exists \eps_0>0\;
    \forall \eps \in (0,\eps_0)\;
    \forall \hat v \in \Sobn^1(\Omega) \colon \;
    \alpha \norm {\hat v}_{\Sob^1(\Omega)}^2
    \le \fra_\eps[\hat v] + \omega \norm {\hat v}^2
  \end{equation*}
  does not hold.
\end{lemma}
\begin{proof}
  Without loss of generality we can assume that $\omega > 0$ in the
  definition of uniform ellipticity.  If $(\fra_\eps)_\eps$ was
  uniformly elliptic, then for any $\hat v\in\dom\mathsf{A}_\eps$
  \begin{align*}
    \norm{\nabla \hat v}_{\Lsqr(\Omega;\dC^d)}^2
    \le \frac 1\alpha
    \Bigl(\fra_\eps[\hat v] + \omega \norm{\hat v}^2\Bigr)
    \le \frac 1\alpha
    \Bigl(\norm{(\sfA_\eps + \ii)\hat v}^2 + \omega \norm{\hat v}^2\Bigr)
    \le \frac {1+\omega}\alpha \norm{(\sfA_\eps + \ii)\hat v}^2
  \end{align*}
  as in the proof of Lemma~\ref{lem:coercive.est}.  In particular, we would have
  \begin{align*}
    \norm{\hat v}_{\Sob^1(\Omega)}^2
    \le \Bigl(\frac{1+\omega}\alpha+1\Bigr)\norm{(\sfA_\eps + \ii)\hat v}^2
  \end{align*}
  and the claim of Lemma~\ref{lem:coercive.est} would follow (with
  another constant).
\end{proof}

%
\section{Absence of norm resolvent convergence}
\label{sec:counterex}
%
The aim of this section is to construct a counterexample to norm
resolvent convergence of the operators $\sfA_\eps$ to the operator
$\sfA$ in the case of resonant potentials and thus to justify that we
can only prove strong resolvent convergence in this setting. This
counterexample relies on the analysis of convergence on the unit disk
$\cB\subset\dR^2$. The model on the disk admits separation of
variables in polar coordinates and the analysis significantly
simplifies. We expect that also for more general domains one can not
hope for norm resolvent convergence of $\sfA_\eps$ to $\sfA$.

In order to construct the counterexample we need to restrict further
the class of resonant potentials. This restriction is clarified in the
following hypothesis.
\begin{hypothesis}
  \label{hypothesis}
  Assume that the resonant potential $V\in \Cci(\clo \dR_+)$ (in the
  sense of Definition~\ref{def:resonant}) is such that the
  self-adjoint one-dimensional Schr\"odinger operator with domain
  $\Sobn^1(\dR_+)\cap \Sob^2(\dR_+)$ acting as
  $\psi \mapsto -\psi'' + V\psi$ in the Hilbert space $\Lsqr(\dR_+)$
  has at least one negative eigenvalue. We denote by $\mu < 0$ the
  lowest eigenvalue of this Schr\"odinger operator and by
  $f_\mu\in \Sobn^1(\dR_+)\cap \Sob^2(\dR_+)$ the corresponding
  real-valued eigenfunction.
\end{hypothesis}

\begin{example}
  Let $V\in \Cci(\clo \dR_+)$ such that $\supp V\subset[0,a]$ with
  $a > 0$ be a non-positive resonant potential in the sense of
  Definition~\ref{def:resonant}; \ie the self-adjoint Schr\"odinger
  operator in $\Lsqr((0,a))$ corresponding to the quadratic form
  $f\mapsto \int_0^a(\abs{f'}^2+V\abs f^2)\dd x$ with domain
  $\{f\in\Sob^1((0,a))\colon f(0) = 0\}$ has eigenvalue zero. Recall
  that there exists a sequence of real numbers
  $\{\alpha_n\}_{n\in\dN}$,
  $1=\alpha_1 < \alpha_2 <\alpha_3 < \dots< \alpha_n <\dots$ such that
  $\alpha_n\to\infty$, for which the multiple $\alpha_n V$ of the
  potential $V$ is resonant for all $n\in\dN$. It remains to note that
  for all $n\in\dN$ sufficiently large the resonant potential
  $\alpha_n V$ necessarily satisfies
  Hypothesis~\ref{hypothesis}. Thus, the family of resonant potentials
  satisfying the above hypothesis is non-void.
\end{example}

The quadratic form of the operator $\sfA_\eps$ in the case of the unit
disk can be written in polar coordinates
\begin{equation*}
  \fra_\eps[u]
  = \int_0^1\int_0^{2\pi}
  \Bigl(\abs{\partial_r u}^2 + \frac{\abs{\partial_\theta u}^2}{r^2}
  + \frac{1}{\eps^2}V\bigg(\frac{1-r}{\eps}\bigg)\abs u^2
  \Bigr) r\dd\theta\dd r,
\end{equation*}
where the form domain remains the Sobolev space $\Sobn^1(\cB)$.
For $m\in\dZ$, consider the quadratic form of the fibre operator:
\begin{align*}
  \fra_{\eps}^{(m)}[f]
  &\defeq \int_0^1\Bigl(\abs{f'(r)}^2
    +\frac{m^2}{r^2}\abs{f(r)}^2
    + \frac{1}{\eps^2}V\bigg(\frac{1-r}{\eps}\bigg)\abs{f(r)}^2
    \Bigr)r\dd r,\\
  \dom\fra_\eps^{(m)}
  & \defeq \big\{ f \in \Lsqr((0,1);r\dd r) \colon
    \fra_\eps^{(m)}[f]<\infty \big\}.
\end{align*}
The symmetric quadratic form $\fra_\eps^{(m)}$ is closed, densely
defined, and lower-semibounded in $L^2((0,1);r\dd r)$ for any $m\in\dZ$. The mentioned
properties of $\fra_\eps^{(m)}$ follow immediately from the
perturbation result~\cite[Chapter VI, Theorem 1.33]{kato:66} and the fact
that this form can be represented as a sum of a bounded quadratic form
\begin{equation*}
  f \mapsto \frac{1}{\eps^{2}}
  \int_0^1 V\bigg(\frac{1-r}{\eps}\bigg) \abs{f(r)}^2r\dd r
 \end{equation*}
 on $\Lsqr((0,1);r\dd r)$ and the quadratic for the fibre operator of
 the Dirichlet Laplacian on the disk, for which these properties are
 well known.  Let us denote by $\sfA_\eps^{(m)}$ the self-adjoint
 fibre operator in $\Lsqr((0,1);r\dd r)$ associated with the quadratic
 form $\fra_\eps^{(m)}$.  Using standard procedure based on separation
 of variables we infer the following unitary equivalence
 \begin{equation}
   \label{eq:unitary_equiv}
   \sfA_\eps \cong \bigoplus_{m\in\dZ} \sfA^{(m)}_\eps.
\end{equation}
In particular, we get as a direct consequence
\begin{equation}
  \label{eq:spec}
  \sigma(\sfA_\eps)
  = \clo{\bigcup_{m\in\dZ} \sigma(\sfA^{(m)}_\eps)}.
\end{equation}
The spectrum of the fibre operator $\sfA_\eps^{(m)}$ is clearly purely
discrete and let us denote by $\lambda_1^{(m)}(\eps)$ the lowest
eigenvalue of $\sfA^{(m)}_\eps$.

The following lemma is essential in the construction of the
counterexample. Its proof is outsourced to Appendix~\ref{app:A}.
\begin{lemma}
  \label{lem:fibres}
  For any $m\in\dZ$, the following properties hold.
  \begin{enumerate}
  \item
    \label{fibres.a}
    $\lambda_1^{(m)}(\cdot)$ is a continuous function.
  \item
    \label{fibres.b}
    $\lim_{\eps \to 0}\lambda_1^{(m)}(\eps) = -\infty$
    for any $V$ satisfying Hypothesis~\ref{hypothesis}.
  \item
    \label{fibres.c}
    $\lambda_1^{(m)}(\eps) \ge 0$ if $\eps
    \ge \frac{1}{\abs m}\sqrt{\norm V_\infty}$.
  \end{enumerate}
\end{lemma}

The next proposition provides a counterexample based on the disk. The
proposed technique can be also used to construct counterexamples for
domains other than the disk.
\begin{proposition}
  \label{prp:counterex}
  For the unit disk and a resonant potential $V$ satisfying
  Hypothesis~\ref{hypothesis}, the family of operators $\sfA_\eps$
  does not converge in norm resolvent sense to the operator $\sfA$.
\end{proposition}
\begin{proof}
  Recall that the Robin Laplacian $\sfA$ is bounded from below. Let us
  choose $\beta <0$ such that $\beta < \inf\sigma(\sfA)$.  By
  Lemma~\ref{lem:fibres}
  we can find $m_1\in\dN$
  and $\eps_1 > 0$ such that $\lambda_1^{(m_1)}(\eps_1) =\beta$.  By
  item~\itemref{fibres.c} of the same lemma we can choose integer
  $m_2 > m_1$ such that $\lambda_1^{(m_2)}(\eps_1) \ge 0$.  Hence, by
  items~\itemref{fibres.a} and~\itemref{fibres.b} of
  Lemma~\ref{lem:fibres} we can find $\eps_2\in(0,\eps_1)$ such that
  $\lambda_1^{(m_2)}(\eps_2) = \beta$. Analogously, we can find
  $m_3 > m_2$ and $\eps_3\in (0,\eps_2)$ such that
  $\lambda_1^{(m_3)}(\eps_3) = \beta$.  Thus, repeating the
  construction, we conclude that there exists sequences of real
  numbers $\eps_1 > \eps_2 > \dots > \eps_k >\dots > 0$ and integers
  $m_1 < m_2 <\dots < m_k < \dots <+\infty$ such that
  $\lambda_1^{(m_k)}(\eps_k) =\beta$ for all $k\in\dN$.  Moreover, it
  follows from Lemma~\ref{lem:fibres}~\itemref{fibres.c} that
  $\eps_k \le \frac{1}{m_k}\sqrt{\norm V_\infty}\to 0$ as $k\to\infty$.

  Suppose for the moment that $\sfA_\eps$ converges to $\sfA$ in the
  norm resolvent sense as $\eps\to0^+$. Then also $\sfA_{\eps_k}$
  converges to $\sfA$ in norm resolvent sense as
  $k\to\infty$. By~\cite[Satz 9.24~(i)]{weidmann:00} we would get that
  the spectrum of the operator $\sfA_{\eps_k}$ must converge to the
  spectrum of the operator $\sfA$ as $k\to\infty$. This consequence of
  the norm resolvent convergence combined with~\eqref{eq:spec}
  contradicts the choice of the sequence $(\eps_k)_{k\in\dN}$, since
  $\beta < \inf\sigma(\sfA)$ is in the spectrum of $\sfA_{\eps_k}$ for
  all $k\in\dN$.
\end{proof}

\subsection*{Acknowledgements}
The first named author (VL) is grateful to Casa Mathematica Oaxaca
(CMO) for providing the possibility to participate in the workshop
``Analytic and Geometric Aspects of Spectral Theory'' in August 2022,
where this project was initiated and to University of Trier for
funding a research stay in February 2023, during which a part of the
results presented in this manuscript were obtained.  The research was
supported by the European Union's Horizon~2020 research and innovation
programme under the Marie Sk\l odowska-Curie grant agreement No 873071.
We are grateful to the anonymous referees for valuable suggestions.

\appendix

%
\section{Proof of Lemma~\ref{lem:fibres}}
\label{app:A}
%

\begin{proof}[Proof of Lemma~\ref{lem:fibres}]
  \itemref{fibres.a}~Let $\eps_0\in(0,\infty)$. It is straightforward
  to see that
  $\lambda_1^{(m)}(\eps) \ge -\frac{4}{\eps_0^2}\|V\|_\infty$ for all
  $\eps\in(\eps_0/2,2\eps_0)$. In other words, the lowest eigenvalue
  $\lambda_1^{(m)}(\eps)$ is uniformly bounded from below for
  $\eps\in(\eps_0/2,2\eps_0)$.  Thus, in view of~\cite[Satz
  9.24]{weidmann:00}, continuity of
  $\eps\mapsto \lambda_1^{(m)}(\eps)$ for any $m\in\dZ$ would
  immediately follow if we show that the operators $\sfA_\eps^{(m)}$
  converge in norm resolvent sense to $\sfA_{\eps_0}^{(m)}$ as
  $\eps \to \eps_0$. To this aim notice that for any
  $\eps\in(\eps_0/2,2\eps_0)$
  \begin{align*}
    \int_0^1 \abs[\Big]{\frac{1}{\eps_0^2}
    V\bigg(\frac{1-r}{\eps_0}\bigg)
    &-\frac{1}{\eps^2} V\bigg(\frac{1-r}{\eps}\bigg)
      } \abs{f(r)}^2 r\dd r\\
    &\le \sup_{r\in(0,1)}
      \abs[\Big]{\frac{1}{\eps^2_0} V\bigg(\frac{1-r}{\eps_0}\bigg)-
      \frac{1}{\eps^2} V\bigg(\frac{1-r}{\eps}\bigg)}
      \int_0^1|f(r)|^2r \dd r\\
    &\le
      \frac{16}{\eps_0^4}\big(\eps_0\norm V_\infty + \norm{V'}_\infty\big)
      \abs{\eps-\eps_0} \int_0^1 \abs{f(r)}^2 r\dd r.
  \end{align*}
  Thus, it follows that
  \begin{equation*}
    \abs[\big]{\fra_\eps^{(m)}[f] - \fra_{\eps_0}^{(m)}[f]}
    \le C\abs{\eps -\eps_0}\int_0^1 \abs{f(r)}^2r\dd r
  \end{equation*}
  for any $\eps\in(\eps_0/2,2\eps_0)$ and all
  $f\in \dom\fra_\eps^{(m)} = \dom\fra_{\eps_0}^{(m)}$ with constant
  $C = C(V,\eps_0) = (16\eps_0^{-4})(\eps_0\norm V_\infty +
  \norm{V'}_\infty) > 0$.  Hence, the norm resolvent convergence of
  $\sfA_\eps^{(m)}$ to $\sfA_{\eps_0}^{(m)}$ as $\eps\to\eps_0$ is a
  consequence of~\cite[Chapter VI, Theorem 3.4]{kato:66}.

  \itemref{fibres.b}~Let the cut-off function $\chi\in \Cci((0,1])$ be
  such that $0\le\chi \le 1$, $\chi(r) = 1$ for $r\in[3/4,1]$, and
  $\chi(r) = 0$ for all $r\in (0,1/2]$.  As a trial function for the
  quadratic form $\fra_\eps^{(m)}$, we use
  \begin{equation*}
    g_\eps(r) := \chi(r)f_\mu\bigg(\frac{1-r}{\eps}\bigg),\qquad r\in(0,1),
  \end{equation*}
  where $f_\mu\in \Sobn^1(\dR_+)\cap \Sob^2(\dR_+)$ satisfies the
  differential equation $-f_\mu'' +V f_\mu = \mu f_\mu$ and where
  $\mu<0$ is as in Hypothesis~\ref{hypothesis}. Since
  $\supp V\subset[0,a]$, we obtain that for $t > a$
  \begin{equation}
    \label{eq:largex}
    f_\mu(t) = C_\mu e^{-t\sqrt{-\mu}}
  \end{equation}
  for some $C_\mu\in\dR\setminus\{0\}$.  By a direct computation, we
  obtain for the square of the weighted $\Lsqr$-norm of $g_\eps$ the
  following asymptotic expansion
  \begin{align}
    \nonumber
    \int_0^1 \abs{g_\eps(r)}^2 r\dd r
    &= \eps\int_0^{1/\eps}
      \chi^2(1-\eps t) \abs{f_\mu(t)}^2(1-\eps t)\dd t\\
    \nonumber
    &= \eps\int_0^{\infty}
      \abs{f_\mu(t)}^2\dd t - \eps\int_{1/\eps}^\infty \abs{f_\mu(t)}^2\dd t
      - \eps^2\int_0^{1/(4\eps)} \abs{f_\mu(t)}^2t\dd t\\
    \nonumber
    &\qquad\qquad
      +\eps\int_{1/(4\eps)}^{1/\eps} ((1-\eps t)
      \chi^2(1-\eps t) -1)\abs{f_\mu(t)}^2\dd t \\
    \label{eq:gL2}
    &= \eps\int_0^\infty \abs{f_\mu(t)}^2\dd t+ o(\eps),
      \qquad\eps \to 0,
  \end{align}
  where in the first step we perform the change of variable
  $r=1-\eps t$, in the second step we decompose the integral term into
  the sum of four integral terms via an identical transform based on
  the properties of $\chi$, and in the last step we used that
  $\abs{(1-\eps t)\chi^2(1-\eps t)-1}\le 1$ for all
  $t\in (1/(4\eps),1/\eps)$ and that
  \begin{equation*}
    \lim_{\eps \to 0}\int_{1/\eps}^\infty \abs{f_\mu(t)}^2\dd t =
    \lim_{\eps \to 0}
    \int_{1/(4\eps)}^{1/\eps} \abs{f_\mu(t)}^2\dd t = 0,
    \qquad
    \int_{0}^\infty \abs{f_\mu(t)}^2t \dd t < \infty,
  \end{equation*}
  where the last integral is finite due to~\eqref{eq:largex}.

  Without loss of generality we may assume in the rest of the argument
  that $\eps < 1/16$.  For the quadratic form $\fra_\eps^{(m)}$ of the
  fibre operator evaluated on the trial function $g_\eps$ we obtain
  using the properties of the cut-off function $\chi$ and the
  substitution $r = 1-t\eps$ that
  \begin{align}
    \nonumber
    \fra_\eps^{(m)}[g_\eps]
    &= \int_0^1\bigg[\bigg(
      \chi'(r) f_\mu\bigg(\frac{1-r}{\eps}\bigg)
      -\frac{1}{\eps} \chi(r) f_\mu'\bigg(\frac{1-r}{\eps}\bigg)\bigg)^2
      +\frac{m^2}{r^2}\chi^2(r)f_\mu^2\bigg(\frac{1-r}{\eps}\bigg)\\
    \nonumber
    &\qquad\qquad
      +\frac{1}{\eps^2}V\bigg(\frac{1-r}{\eps}\bigg)
      \chi^2(r)f_\mu^2\bigg(\frac{1-r}{\eps}\bigg)\bigg]r\dd r\\
    \nonumber
    &\le \eps\int_0^{1/\eps}
      \bigg[\Big(\chi'(1-t\eps) f_\mu(t)
      -\frac{1}{\eps} \chi(1-\eps t) f_\mu'(t)\Big)^2
      +4m^2\chi^2(1-\eps t)\abs{f_\mu(t)}^2\\
    \label{eq:a1}
    &\qquad\qquad+\frac{1}{\eps^2}V(t)\chi^2(1-\eps t)\abs{f_\mu(t)}^2
      \bigg](1-\eps t)\dd t = I(\eps) + J(\eps),
  \end{align}
  where the terms $I(\eps)$ and $J(\eps)$ are defined by
  \begin{align*}
    I(\eps)
    &:= \eps\int_0^{1/(\sqrt{\eps})}
      \bigg[\frac{1}{\eps^2} \abs{f_\mu'(t)}^2
      + 4m^2 \abs{f_\mu(t)}^2
      +\frac{1}{\eps^2} V(t)\abs{f_\mu(t)}^2
      \bigg](1-\eps t)\dd t \quad\text{and}\\
    J(\eps)
    &:= \eps \int_{1/(\sqrt{\eps})}^{1/\eps}
      \bigg[\Big(\chi'(1-t\eps) f_\mu(t)
      -\frac{1}{\eps} \chi(1-\eps t) f_\mu'(t)\Big)^2
      +4m^2\chi^2(1-\eps t)\abs{f_\mu(t)}^2\\
    &\qquad\qquad
      +\frac{1}{\eps^2}V(t)\chi^2(1-\eps t)\abs{f_\mu(t)}^2
      \bigg](1-\eps t)\dd t.
  \end{align*}
  Using that
  \begin{equation*}
    \int_0^{1/(\sqrt{\eps})} \bigg[|f_\mu'(t)|^2
    +V(t)\abs{f_\mu(t)}^2\bigg]\dd t \to
    \int_0^\infty \bigg[|f_\mu'(t)|^2
    +V(t)\abs{f_\mu(t)}^2\bigg]\dd t = \mu\int_0^\infty \abs{f_\mu(t)}^2\dd t
  \end{equation*}
  as $\eps \to 0$, we obtain
  \begin{align}
    \nonumber
    I(\eps)
    &= \frac{1}{\eps}\int_0^{1/(\sqrt{\eps})}
      \bigg[
      \abs{f_\mu'(t)}^2
      + 4m^2\eps^2 \abs{f_\mu(t)}^2
      + V(t)\abs{f_\mu(t)}^2
      \bigg]\dd t\\
    \nonumber
    &\qquad\qquad -
      \int_0^{1/(\sqrt{\eps})}
      \bigg[
      \abs{f_\mu'(t)}^2
      + 4m^2 \eps^2 \abs{f_\mu(t)}^2
      + V(t)\abs{f_\mu(t)}^2
      \bigg] t\dd t\\
    \label{eq:I}
    &= \frac{\mu}{\eps}\int_0^\infty \abs{f_\mu(t)}^2 \dd t
      + o(\eps^{-1}),
      \qquad\eps \to 0.
  \end{align}
  Moreover, we conclude applying the inequality $(a+b)^2\le 2a^2+2b^2$
  (valid for any $a,b > 0$) and the properties of $\chi$ that
  \begin{align}
    \nonumber
    \abs{J(\eps)}
    &\le \eps\int_{1/\sqrt{\eps}}^{1/\eps}
      \bigg[
      2\norm{\chi'}_\infty^2 \abs{f_\mu(t)}^2
      + \frac{2}{\eps^2}\abs{f_\mu'(t)}^2
      + 4m^2 \abs{f_\mu(t)}^2
      + \frac{1}{\eps^2}\norm V_\infty \abs{f_\mu(t)}^2
      \bigg]\dd t\\
    \label{eq:J}
    &= o(\eps^{-1}).
  \end{align}
  Plugging~\eqref{eq:I} and~\eqref{eq:J} into~\eqref{eq:a1} we end up
  with the asymptotic expansion
  \begin{equation}
    \label{eq:a2}
    \fra_\eps^{(m)}[g_\eps]
    = \frac{\mu}{\eps}\int_0^\infty \abs{f_\mu(t)}^2\dd t+ o(\eps^{-1}),
    \qquad\eps \to 0.
  \end{equation}
  Finally, combining~\eqref{eq:gL2} and~\eqref{eq:a2} with the min-max
  principle we arrive at
  \begin{equation*}
    \lambda_1^{(m)}(\eps)
    \le \frac{\fra_\eps^{(m)}[g_\eps]}{\int_0^1\abs{g_\eps(r)}^2r\dd r}
    = \frac{\mu}{\eps^2} + o(\eps^{-2}),
    \qquad \eps\to 0.
  \end{equation*}
  The claim then follows from the fact that $\mu < 0$.

  \itemref{fibres.c}~The statement is a consequence of the
  representation of the quadratic form $\fra_{\eps}^{(m)}$ and the
  fact that under the assumption
  $\eps\ge (1/\abs m) \sqrt{\norm V_\infty}$ the function on $(0,1)$
  given by
   \begin{equation*}
    r \mapsto
    \frac{m^2}{r^2} + \frac{1}{\eps^2}V\bigg(\frac{1-r}{\eps}\bigg)
  \end{equation*}
  is non-negative.
\end{proof}

%

\newcommand{\etalchar}[1]{$^{#1}$}
\providecommand{\bysame}{\leavevmode\hbox to3em{\hrulefill}\thinspace}
\providecommand{\MR}{\relax\ifhmode\unskip\space\fi MR }
\providecommand{\MRhref}[2]{%
  \href{http://www.ams.org/mathscinet-getitem?mr=#1}{#2}
}
\providecommand{\href}[2]{#2}

\end{document}